\documentclass[12pt]{article}

\usepackage{color}
\usepackage{amssymb}
\usepackage{amsmath}
\usepackage{cite}
\usepackage[arrow, matrix, curve]{xy}

\begin{document}

\renewcommand{\citeleft}{{\rm [}}
\renewcommand{\citeright}{{\rm ]}}
\renewcommand{\citepunct}{{\rm,\ }}
\renewcommand{\citemid}{{\rm,\ }}

\newcounter{abschnitt}
\newtheorem{satz}{Theorem}
\newtheorem{theorem}{Theorem}[abschnitt]
\newtheorem{koro}[theorem]{Corollary}
\newtheorem{prop}[theorem]{Proposition}
\newtheorem{lem}[theorem]{Lemma}
\newtheorem{expls}[theorem]{Examples}
\newtheorem{expl}[theorem]{Example}

\renewenvironment{quote}{\list{}{\leftmargin=0.62in\rightmargin=0.62in}\item[]}{\endlist}

\newcounter{saveeqn}
\newcommand{\alpheqn}{\setcounter{saveeqn}{\value{abschnitt}}
\renewcommand{\theequation}{\mbox{\arabic{saveeqn}.\arabic{equation}}}}
\newcommand{\reseteqn}{\setcounter{equation}{0}
\renewcommand{\theequation}{\arabic{equation}}}

\hyphenpenalty=9000

\sloppy

\phantom{a}

\vspace{-1.7cm}

\begin{center}
\begin{Large} {\bf Minkowski Valuations and Generalized Valuations} \\[0.7cm] \end{Large}

\begin{large} Franz E. Schuster and Thomas Wannerer \end{large}
\end{center}

\vspace{-0.8cm}

\begin{quote}
\footnotesize{ \vskip 1cm \noindent {\bf Abstract.}
A convolution representation of continuous translation invariant and
$\mathrm{SO}(n)$ equivariant Minkowski valuations is established. This is based on a new classification
of translation invariant generalized spherical valuations. As applications, Crofton and kinematic formulas for
Minkowski valuations are obtained. }
\end{quote}

\vspace{0.6cm}

\centerline{\large{\bf{ \setcounter{abschnitt}{1}
\arabic{abschnitt}. Introduction}}}

\alpheqn

\vspace{0.6cm}

A valuation on convex bodies (non-empty compact convex sets) is a \emph{finitely additive} function.
More precisely, let $\mathcal{K}^n$ denote the space of convex bodies in $\mathbb{R}^n$ endowed with the Hausdorff metric.
A map $\phi: \mathcal{K}^n \rightarrow \mathcal{A}$ with values in an Abelian semigroup $\mathcal{A}$ is a {\it valuation} if
\[\phi(K) + \phi(L) = \phi(K \cup L) + \phi(K \cap L)   \]
whenever $K \cup L$ is convex. As a generalization of the notion of measure and as the crucial ingredient in Dehn's solution of Hilbert's third problem,
scalar valuations (where $\mathcal{A} = \mathbb{R}$ or $\mathbb{C}$) have long played a central role in convex and discrete geometry
(see \textbf{\cite{Klain:Rota}} or \textbf{\cite[\textnormal{Chapter 6}]{schneider93}}). The most famous classical result on valuations is a celebrated characterization of rigid motion invariant valuations by Hadwiger \textbf{\cite{hadwiger51}} (which was slightly improved later by Klain \textbf{\cite{klain95}}).

\begin{satz}[\!\!\cite{hadwiger51, klain95}] \label{hadwiger} The intrinsic volumes $V_0, V_1, \ldots, V_n$ form a basis
of the vector space of all continuous scalar valuations on $\mathcal{K}^n$ which are translation and $\mathrm{SO}(n)$ invariant.
\end{satz}

Hadwiger's characterization theorem had a transformative effect on integral geometry. It not only allowed
for almost effortless proofs of the principal and more general kinematic formulas (see, e.g., \textbf{\cite{Klain:Rota}}) but also
made the importance of precise descriptions of classes of invariant valuations evident. Still to this day, Theorem \ref{hadwiger} often serves as
a starting point for the classification of invariant scalar valuations (see, e.g., \textbf{\cite{Alesker01, AleskerFaifman, Bernig09, habparap14, centro}}) and, more general, equivariant tensor valued valuations, where $\mathcal{A} = \mathrm{Sym}^k\mathbb{R}^n$ (see \textbf{\cite{Alesker99, ABS2011, wannerer14, hugschnschu07}}).
These results in turn were critical for the tremendous progress in integral geometry of recent years \linebreak (see \textbf{\cite{Alesker03, bernigfu10, fu06, bernighug2015, hugschnschu07, wannerer13}} and the references therein).

\vspace{0.2cm}

In 1974 Schneider \textbf{\cite{schneider74, schneider74b}} first investigated valuations, where $\mathcal{A} = \mathcal{K}^n$ and addition on $\mathcal{K}^n$ is the usual Minkowski addition.
In a more recent influential article, Ludwig \textbf{\cite{Ludwig:Minkowski}} coined the name \emph{Minkowski valuations} for such maps and started a line of research
concerned with the classification of Minkowski valuations intertwining \emph{linear} transformations, see \textbf{\cite{abardia12, abardber11, haberl11, ludwig02, Ludwig10a, SchuWan11, wannerer10}}.

\pagebreak

The recent results on Minkowski valuations which are equi- or contravariant with respect to linear transformations show that they often form convex cones generated by \emph{finitely} many valuations,
such as the projection or difference body operators. In contrast to this, the cone of translation invariant and $\mathrm{SO}(n)$ equivariant Minkowski valuations is \emph{infinite dimensional}.
This is one reason why no full analogue of Theorem \ref{hadwiger} for Minkowski valuations has been obtained yet, except for dimension
$n = 2$, where Schneider \textbf{\cite{schneider74b}} already established such a result. We therefore assume throughout that $n \geq 3$.

About a decade ago, Kiderlen \textbf{\cite{kiderlen05}} and the first author \textbf{\cite{Schu06a}} were the first to obtain convolution representations
of translation invariant and $\mathrm{SO}(n)$ equivariant \emph{continuous} Minkowski valuations. However, their results were limited to the case of valuations of
degree $1$ and $n - 1$, respectively, where a map $\Phi$ from $\mathcal{K}^n$ to $\mathcal{K}^n$ (or
$\mathbb{R}$) is said to have {\it degree} $i$ if $\Phi(\lambda K) = \lambda^i\Phi K$ for $K \in \mathcal{K}^n$ and $\lambda >
0$. The convolution of functions and measures on $\mathbb{S}^{n-1}$ used in \textbf{\cite{kiderlen05}} and \textbf{\cite{Schu06a}} is induced from the group $\mathrm{SO}(n)$ by identifying
$\mathbb{S}^{n-1}$ with the homogeneous space $\mathrm{SO}(n)/\mathrm{SO}(n - 1)$ (see Section 2 for details).

Under additional \emph{smoothness} assumptions, the first author in \textbf{\cite{Schu09}} and jointly with the second author in \textbf{\cite{SchuWan13}} extended the
results from \textbf{\cite{kiderlen05}} and \textbf{\cite{Schu06a}} to the remaining (non-trivial) degrees $i \in \{2, \ldots, n - 2\}$ when the Minkowski valuations are \emph{even}. (McMullen \textbf{\cite{McMullen77}} showed that only integer degrees $0 \leq i \leq n$ can occur.) However, the techniques employed in \textbf{\cite{Schu09}} or \textbf{\cite{SchuWan13}} were not suited to describe  merely continuous Minkowski valuations, which is the goal since the 1970s.

\vspace{0.2cm}

In this article we establish a precise description of all \emph{continuous} translation invariant and
$\mathrm{SO}(n)$ equivariant Minkowski valuations without any further assumptions on the parity or degree of the valuations.
As we explain in Section~5, our main theorem generalizes and implies all previously obtained
convolution representations of Minkowski valuations intertwining rigid motions. In order to state our result, recall
that a convex body $K \in \mathcal{K}^n$ is uniquely determined by its support function $h_K(u)=\max \{u \cdot x: x \in K\}$ for $u \in
\mathbb{S}^{n-1}$ and let  $\mathcal{M}_{\mathrm{o}}(\mathbb{S}^{n-1})$ and $C_{\mathrm{o}}(\mathbb{S}^{n-1})$ denote the spaces of signed Borel measures
and continuous functions on $\mathbb{S}^{n-1}$, respectively, having their center of mass at the origin.

\begin{satz} \label{main1} If $\Phi: \mathcal{K}^n \rightarrow \mathcal{K}^n$ is a continuous Minkowski valuation which
is translation invariant and $\mathrm{SO}(n)$ equivariant, then there exist uniquely determined $c_0, c_n \geq 0$,
$\mathrm{SO}(n-1)$ invariant $\mu_i \in \mathcal{M}_{\mathrm{o}}(\mathbb{S}^{n-1})$, $1 \leq i \leq n - 2$, and an $\mathrm{SO}(n-1)$ invariant
$f_{n-1} \in C_{\mathrm{o}}(\mathbb{S}^{n-1})$ such that
\begin{equation} \label{mainform}
h_{\Phi K} = c_0 + \sum_{i=1}^{n-2} S_i(K,\cdot) \ast \mu_i + S_{n-1}(K,\cdot) \ast f_{n-1} + c_n V_n(K)
\end{equation}
for every $K \in \mathcal{K}^n$.
\end{satz}

\pagebreak

The Borel measures $S_i(K,\cdot)$, $1 \leq i \leq n - 1$, on $\mathbb{S}^{n-1}$ are Aleksandrov's \emph{area measures} (see, e.g., \textbf{\cite{schneider93}})
associated with $K \in \mathcal{K}^n$. If $K$ is sufficiently smooth and has positive curvature, then each $S_i(K,\cdot)$ is absolutely continuous with
respect to spherical Lebesgue measure and its density is (up to a constant) given by the $i$th elementary symmetric function of the
principal radii of curvature of $K$.

The equality in (\ref{mainform}) has to be understood in the sense of measures, where we identify $h \in C(\mathbb{S}^{n-1})$ with the absolutely continuous measure with density $h$. \linebreak
For $n \leq 4$, we show in Section 5 that if $\Phi$ has degree $1$ or $2$, then the measures $\mu_1$ or $\mu_2$, respectively,
are in fact absolutely continuous with a density in $L^2(\mathbb{S}^{n-1})$. However, this is no longer true in general when $n > 4$.

\vspace{0.2cm}

The proof of Theorem \ref{main1} is based on new techniques involving translation invariant \emph{generalized} valuations which
were only recently introduced by Alesker and Faifman \textbf{\cite{AleskerFaifman}} (see also \textbf{\cite{bernigfaifman15}}).
Generalized valuations are related to smooth valuations in the same way that distributions are related to smooth functions.
More precisely, let $\mathbf{Val}^{\infty}_i$, $0 \leq i \leq n$, denote the space of \emph{smooth} translation invariant valuations of degree $i$
endowed with the G{\aa}rding topology which makes it a Fr\'echet space (see Section 3 for more information). The space $\mathbf{Val}^{-\infty}_i$ of (translation invariant) \emph{generalized valuations}
of degree $i$ is defined as the topological dual
\begin{equation} \label{defgenval}
\mathbf{Val}^{-\infty}_i := (\mathbf{Val}_{n-i}^{\infty})^*
\end{equation}
endowed with the weak topology.

As part of his far reaching reconceptualization of integral geometry, Alesker \textbf{\cite{Alesker04a}} discovered
a continuous \emph{non-degenerate} bilinear pairing
\[\langle\,\cdot\,,\cdot\,\rangle: \mathbf{Val}^{\infty}_i \times \mathbf{Val}^{\infty}_{n-i} \rightarrow \mathbb{R} \]
for $0 \leq i \leq n$ (see also Section 3). The induced \emph{Poincar\'e duality map}
\begin{equation*}
\mathrm{pd}: \mathbf{Val}^{\infty}_i \rightarrow  (\mathbf{Val}_{n-i}^{\infty})^* = \mathbf{Val}^{-\infty}_i
\end{equation*}
is therefore continuous, injective and has dense image with respect to the weak topology. This was the motivation for definition (\ref{defgenval}) and shows
that $\mathbf{Val}^{-\infty}_i$ can be seen as a completion of $\mathbf{Val}^{\infty}_i$ with respect to the weak topology.
Alesker \textbf{\cite{alesker10}} also proved that the Poincar\'e duality map admits a unique continuous extension to the space $\mathbf{Val}_i$
of \emph{continuous} translation invariant valuations of degree $i$. Thus, just like smooth and continuous functions or, more general, signed Borel measures can be identified with
subclasses of distributions (compare Section 2), we can use the Poincar\'e duality map in the following to \emph{identify} the spaces
$\mathbf{Val}_i^{\infty}$ or $\mathbf{Val}_i$, respectively, with certain dense subspaces of $\mathbf{Val}^{-\infty}_i$.

It was first observed in \textbf{\cite{Schu09}} that a translation invariant and $\mathrm{SO}(n)$ equivariant continuous Minkowski valuation $\Phi$
is uniquely determined by a scalar valuation $\varphi \in \mathbf{Val}_i^{\mathrm{SO}(n - 1)}$, the subspace of $\mathrm{SO}(n-1)$ invariant valuations in $\mathbf{Val}_i$.
In turn, valuations in $\mathbf{Val}_i^{\mathrm{SO}(n - 1)}$ are \emph{spherical}.

\pagebreak

Spherical (generalized) valuations correspond
to \emph{spherical representations} of $\mathrm{SO}(n)$ (see Section 3 for details). Let $\mathbf{Val}_i^{\infty,\mathrm{sph}}$
and $\mathbf{Val}_i^{-\infty,\mathrm{sph}}$ denote the subspaces of smooth and generalized spherical valuations, respectively, and let $C_{\mathrm{o}}^{-\infty}(\mathbb{S}^{n-1})$
denote the space of distributions on $\mathbb{S}^{n-1}$ which vanish on restrictions of linear functions to $\mathbb{S}^{n-1}$.
Our second main result, which is critical for the proof of Theorem~\ref{main1} but also of independent interest, is the following classification
of (generalized) spherical valuations.

\begin{satz} \label{main2} Suppose that $1 \leq i \leq n - 1$.
\begin{enumerate}
\item[(a)] The map $\mathrm{E}_i: C_{\mathrm{o}}^{\infty}(\mathbb{S}^{n-1}) \rightarrow \mathbf{Val}_i^{\infty,\mathrm{sph}}$, defined by
\[(\mathrm{E}_if)(K) = \int_{\mathbb{S}^{n-1}} f(u)\,dS_i(K,u),  \]
is an $\mathrm{SO}(n)$ equivariant isomorphism of topological vector spaces which admits a unique extension by continuity in the weak topologies to an isomorphism
\[\widetilde{\mathrm{E}}_i: C_{\mathrm{o}}^{-\infty}(\mathbb{S}^{n-1}) \rightarrow \mathbf{Val}_i^{-\infty,\mathrm{sph}}.   \]
\item[(b)] The space $\mathbf{Val}_i^{\mathrm{SO}(n - 1)}$ is contained in $\widetilde{\mathrm{E}}_i(\mathcal{M}_{\mathrm{o}}(\mathbb{S}^{n-1}))$
if $i \leq n - 2$ and in $\widetilde{\mathrm{E}}_i(C_{\mathrm{o}}(\mathbb{S}^{n-1}))$ if $i = n - 1$.
\end{enumerate}
\end{satz}

Theorem \ref{main2} (a) for $i = 1$ was recently proved by Alesker, see \textbf{\cite[\textnormal{Appendix}]{BPSW2014}}.
Theorem \ref{main2} (b) for $i = n - 1$ follows from a classical result of McMullen \textbf{\cite{McMullen80}}.

\vspace{0.2cm}

Characterizations of Minkowski valuations, in particular, earlier versions of Theorem \ref{main1}, have had far reaching implications for isoperimetric type inequalities
(see, e.g., \textbf{\cite{abardber11, BPSW2014, habschu09, LYZ2000a, LYZ2010a, LYZ2010b, SchuWan11}}). Motivated by a recent important Crofton type formula for the identity map of Goodey and Weil \textbf{\cite{goodeyweil2}},
we show in the final section of this paper how Theorem \ref{main1} can be applied to obtain a general Crofton formula for continuous Minkowski valuations which
generalizes the result from \textbf{\cite{goodeyweil2}} and an earlier result of this type from \textbf{\cite{SchuWan13}}. Using our new Crofton formula and Hadwiger's general integral geometric theorem, a consequence of Theorem \ref{hadwiger} (cf.\ \textbf{\cite[\textnormal{p.\ 173}]{schneiderweil}}), we can then also state a kinematic formula for Minkowski valuations.

\vspace{1cm}

\centerline{\large{\bf{ \setcounter{abschnitt}{2}
\arabic{abschnitt}. Preliminaries}}}

\reseteqn \alpheqn \setcounter{theorem}{0}

\vspace{0.6cm}

In this section we first recall basic notions from Riemannian geometry and compute several quantities
in cylindrical coordinates on $\mathbb{S}^{n-1}$ which will be needed in the proof of Theorem \ref{main2} (b).
Next, we collect background material from representation theory and harmonic analysis, in particular, about the convolution of functions
and measures on $\mathbb{S}^{n-1}$ and its relation to the theory of spherical harmonics. We also
recall some well known facts about distributions on $\mathbb{S}^{n-1}$ and the definition of Berg's
functions used in the solution of the classical Christoffel problem.

\pagebreak

Although in this article we are mainly concerned with the Euclidean unit sphere $\mathbb{S}^{n-1}$ in $\mathbb{R}^n$ and the Lie group $\mathrm{SO}(n)$ of proper
rotations of $\mathbb{R}^n$, let us first consider a general smooth manifold $M$. We denote by $C^{\infty}(M)$ the space of all smooth functions on $M$
equipped with the Fr\'echet space topology of uniform convergence of each finite number of derivatives on each compact subset of $M$.
For a Banach space $X$, the Fr\'echet space $C^{\infty}(M,X)$ of all infinitely differentiable functions on $M$ with values in $X$ is defined similarly.

If in addition $M$ is compact and endowed with a Riemannian metric, then the $C^k$ norm $\|f\|_{C^k}$ of a function $f \in C^k(M)$ is defined by
(see, e.g., \textbf{\cite[\textnormal{p.\ 301}]{petersen}})
\begin{equation} \label{defcknorm}
\|f\|_{C^k} = \sum_{j=0}^k \max \limits_{M} |\nabla^j f|,
\end{equation}
where $\nabla$ denotes the covariant derivative with respect to the given Riemannian metric and $|\nabla^j f|$ is the (Euclidean) norm of
the tensor field $\nabla^j f$.

Among other quantities, we compute in the following example the $C^2$ norm of an $\mathrm{SO}(n - 1)$ invariant function on $\mathbb{S}^{n-1}$ more explicitly.
This will be useful later on for the proof of Theorem \ref{main2} (b).

\vspace{0.3cm}

\begin{expl} \label{exp1} \end{expl}

\vspace{-0.2cm}

\noindent In this article we use $\bar{e} \in \mathbb{S}^{n-1}$ to denote an arbitrary but fixed point (the pole) of the sphere and we write $\mathrm{SO}(n-1)$ for the stabilizer
in $\mathrm{SO}(n)$ of $\bar{e}$. Clearly, every $u \in \mathbb{S}^{n-1}\backslash \{-\bar{e},\bar{e}\}$ can be written uniquely in the form
\begin{equation} \label{cylindcoord}
u = t\bar{e} + \sqrt{1 - t^2} v
\end{equation}
for some $t \in (-1,1)$ and $v \in \mathbb{S}^{n-2}_{\bar{e}} = \{w \in \mathbb{S}^{n-1}: \bar{e} \cdot w = 0\}$.
In the cylindrical coordinates (\ref{cylindcoord}), the (standard) metric tensor $\varrho$ on $\mathbb{S}^{n-1}$ is given by
\begin{equation} \label{metrictens}
\varrho = \frac{1}{1-t^2}\,dt\otimes dt + (1 - t^2)\,dv \otimes dv,
\end{equation}
where $dv \otimes dv$ is the metric tensor on $\mathbb{S}^{n-2}_{\bar{e}}$.

Let $\Delta_{\mathbb{S}}$ denote the Laplacian (or Laplace--Beltrami
operator) on $\mathbb{S}^{n-1}$ and recall that, for $f, g \in C^2(\mathbb{S}^{n-1})$, we have
\begin{equation}
\int_{\mathbb{S}^{n-1}} f(u)\,\Delta_{\mathbb{S}}g(u)\,du = \int_{\mathbb{S}^{n-1}} g(u)\,\Delta_{\mathbb{S}} f(u)\,du.
\end{equation}
Using (\ref{metrictens}), one can easily obtain the following expression for the Laplacian in cylindrical coordinates (cf.\ \textbf{\cite[\textnormal{Proposition 2.6}]{cberg}})
\begin{equation} \label{lapcylin}
\Delta_{\mathbb{S}} = \frac{1}{\sqrt{|\varrho|}}\partial_i \left (\sqrt{|\varrho|} \varrho^{ij}\partial_j \right ) = (1-t^2)\frac{\partial^2}{\partial t^2} - ( n - 1)\,t\frac{\partial}{\partial t} + \frac{1}{1-t^2}\,\bar{\Delta}_{\mathbb{S}},
\end{equation}
where $\bar{\Delta}_{\mathbb{S}}$ denotes the Laplacian on $\mathbb{S}^{n-2}_{\bar{e}}$.

\pagebreak

Now let $f \in C^2(\mathbb{S}^{n-1})$ be $\mathrm{SO}(n - 1)$ invariant, that is, in the cylindrical coordinates (\ref{cylindcoord}) the function $f$ depends only on $t$.
Then, by (\ref{lapcylin}), we have
\begin{equation} \label{lapcylininv}
\Delta_{\mathbb{S}}f = (1-t^2)\frac{\partial^2f}{\partial t^2} - ( n - 1)\,t\frac{\partial f}{\partial t}.
\end{equation}
Moreover, a straightforward computation, using again (\ref{metrictens}), yields
\begin{equation} \label{normnabf}
|\nabla f |^2 = (1 - t^2)\left (\frac{\partial f}{\partial t} \right )^2
\end{equation}
and
\begin{equation} \label{normnab2f}
|\nabla^2 f|^2 = (n - 2)\left (t\,\frac{\partial f}{\partial t} \right )^2 + \left ( (1-t^2)\,\frac{\partial^2 f}{\partial t^2} -t\,\frac{\partial f}{\partial t} \right )^2.
\end{equation}

\vspace{0.4cm}

We turn now to representations of Lie groups. First recall that since $\mathrm{SO}(n)$ is compact all its irreducible representations are finite
dimensional and that the equivalence classes of irreducible complex
representations of $\mathrm{SO}(n)$ are uniquely determined by their highest
weights (see, e.g., \textbf{\cite{broecker_tomdieck}}). These highest weights, in turn, can be indexed by $\lfloor n/2 \rfloor$-tuples of integers
$(\lambda_1,\lambda_2,\ldots,\lambda_{\lfloor n/2 \rfloor})$ such
that
\begin{equation} \label{heiwei}
\left \{\begin{array}{ll} \lambda_1 \geq \lambda_2 \geq \cdots \geq \lambda_{\lfloor n/2 \rfloor} \geq 0 & \quad \mbox{for odd }n, \\
\lambda_1 \geq \lambda_2 \geq \cdots \geq \lambda_{n/2-1} \geq
|\lambda_{n/2}| & \quad \mbox{for even }n.
\end{array} \right .
\end{equation}

A notion of particular importance for our purposes is that of smooth vectors of an infinite-dimensional representation of a Lie group.

\vspace{0.4cm}

\noindent {\bf Definition} \emph{Let $\rho$ be a continuous representation of a Lie group $G$ on a Banach space $X$.
An element $x \in X$ is called a \emph{smooth vector} if the map $z_x: G \rightarrow X$, defined by $z_x(\vartheta)=\rho(\vartheta)x$, is infinitely differentiable.
The subspace of all smooth vectors in $X$ is denoted by $X^{\infty}$.}

\vspace{0.4cm}

It is well known (cf.\ \textbf{\cite[\textnormal{Section 4.4}]{warner72}}) that the subspace $X^{\infty}$ is a $G$ invariant and dense
subset of $X$. Moreover, the map $X^{\infty} \rightarrow C^{\infty}(G,X)$, given by $x \mapsto z_x$, leads to an identification of $X^{\infty}$
with a \emph{closed} subspace of $C^{\infty}(G,X)$. Hence, we can endow $X^{\infty}$ with the relative topology induced by $C^{\infty}(G,X)$.
This topology on $X^{\infty}$ is called the \emph{G{\aa}rding topology} and turns $X^{\infty}$ into a Fr\'echet space. An important property of the G{\aa}rding topology
on $X^{\infty}$ is that the restriction of the representation of $G$ to $X^{\infty}$ is \emph{continuous}.

In the following lemma, we state two more basic facts about smooth vectors which we will use frequently.

\pagebreak

\begin{lem} \label{smoothlem17} Let $G$ be a Lie Group.
\begin{enumerate}
\item[(a)] If $\rho$ and $\tau$ are continuous representations of $G$ on Banach spaces $X$ and $Y$ and $T: X \rightarrow Y$ is
a continuous and $G$ equivariant linear map, then $T(X^{\infty}) \subseteq Y^{\infty}$ and the restriction $T: X^{\infty} \rightarrow Y^{\infty}$ is
continuous.
\item[(b)] If $H$ is a closed subgroup of $G$, then the smooth vectors of the left regular representation of $G$ on $C(G/H)$ are precisely the smooth
functions on $G/H$, that is, $(C(G/H))^{\infty} = C^{\infty}(G/H)$.
\end{enumerate}
\end{lem}

In this article, we are specifically interested in spherical representations of $\mathrm{SO}(n)$ with respect to $\mathrm{SO}(n - 1)$.

\vspace{0.3cm}

\noindent {\bf Definition} \emph{Let $G$ be a compact Lie group and $H$ a closed subgroup of $G$.
A representation $\rho$ of $G$ on a vector space $X$ is called \emph{spherical} with respect to $H$ if there exists
an $H$-invariant non-zero $x \in X$, that is, $\rho(\vartheta)x = x$ for every $\vartheta \in H$.}

\vspace{0.3cm}

For the following two important facts about spherical representations (see, e.g., \textbf{\cite[\textnormal{p.\ 17}]{takeuchi}}), we consider the left
regular representation of $G$ on the Hilbert space $L^2(G/H)$ of square-integrable functions on the homogeneous space $G/H$.

\begin{theorem} \label{thmspher} Let $G$ be a compact Lie group and $H$ a closed subgroup of $G$.
\begin{enumerate}
\item[(i)] Every subrepresentation of $L^2(G/H)$ is spherical with respect to $H$.
\item[(ii)] Every irreducible representation of $G$ which is spherical with respect to $H$ is isomorphic to a subrepresentation of $L^2(G/H)$.
\end{enumerate}
\end{theorem}

\begin{expl} \label{exp2} \end{expl}

\vspace{-0.2cm}

\noindent The decomposition of $L^2(\mathbb{S}^{n-1})$ into an \emph{orthogonal} sum of $\mathrm{SO}(n)$ irreducible subspaces is given by
\begin{equation*} \label{decompsn1}
L^2(\mathbb{S}^{n-1}) = \bigoplus_{k \in \mathbb{N}} \mathcal{H}_k^n,
\end{equation*}
where $\mathcal{H}_k^n$ is the space of spherical harmonics of
dimension $n$ and degree $k$. It is well known that the highest weights corresponding to the spaces $\mathcal{H}_k^n$ are the $\lfloor n/2 \rfloor$-tuples
$(k,0,\ldots,0)$. Since $\mathbb{S}^{n-1}$ is diffeomorphic to the homogeneous space $\mathrm{SO}(n)/\mathrm{SO}(n-1)$, it follows from Theorem \ref{thmspher} that
every irreducible representation of $\mathrm{SO}(n)$ which is spherical with respect to $\mathrm{SO}(n - 1)$ is isomorphic to one of the spaces
$\mathcal{H}_k^n$ and, thus, their highest weights are of the form $(k,0,\ldots,0)$, $k \in \mathbb{N}$.

\vspace{0.4cm}

For the discussion of Theorem \ref{main1} and its applications to integral geometry, we need some more background from the theory of spherical harmonics (see, e.g., \textbf{\cite{groemer96}}).
Let $N(n,k)$ denote the dimension of the space $\mathcal{H}_k^n$ and recall that
\begin{equation} \label{nnk}
N(n,k) = \frac{n + 2k - 2}{n + k - 2} {n + k - 2 \choose n - 2} = \mathrm{O}(k^{n - 2}) \mbox{ as } k \rightarrow \infty.
\end{equation}

\pagebreak

Using $\pi_k: L^2(\mathbb{S}^{n-1}) \rightarrow \mathcal{H}_k^n$ to denote the orthogonal projection, we can write
\begin{equation} \label{fourierexp}
f \sim \sum_{k=0}^{\infty} \pi_k f
\end{equation}
for the (condensed) Fourier expansion of $f \in L^2(\mathbb{S}^{n-1})$. Recall that the Fourier series in (\ref{fourierexp}) converges to $f$ in the $L^2$ norm.

In the theory of spherical harmonics, a function or measure on $\mathbb{S}^{n-1}$ which is $\mathrm{SO}(n-1)$ invariant is often called \emph{zonal}.
The subspace of zonal functions in $\mathcal{H}_k^n$ is $1$-dimensional for every $k \in \mathbb{N}$ and spanned by the function
$u \mapsto P_k^n(u \cdot \bar{e})$, where $P_k^n \in C([-1,1])$ denotes the \emph{Legendre polynomial}  of dimension $n$ and degree~$k$.
Since the spaces $\mathcal{H}_k^n$ are orthogonal, it is not difficult to show that any zonal function $f \in L^2(\mathbb{S}^{n-1})$
admits a series expansion of the form
\begin{equation} \label{expzonal}
f \sim \sum_{k=0}^{\infty} \frac{N(n,k)}{\omega_n}\, a_k^n[f\,]\,P_k^n(\,\,.\cdot \bar{e}),
\end{equation}
where $\omega_n$ denotes the surface area of the $n$-dimensional Euclidean unit ball and
\begin{equation} \label{multleg}
a_k^n[f\,] = \omega_{n-1} \int_{-1}^1 f(t)\,P_k^n(t)\,(1-t^2)^{\frac{n-3}{2}}\,dt.
\end{equation}
Here, we have used again the cylindrical coordinates (\ref{cylindcoord}) to identify the zonal function $f$ with a function on $[-1,1]$.

Now we turn to (formal) Fourier expansions of measures and, more general, distributions on $\mathbb{S}^{n-1}$.
To this end, first recall that a \emph{distribution} on $\mathbb{S}^{n-1}$ is a continuous linear
functional on $C^{\infty}(\mathbb{S}^{n-1})$. Since $\mathbb{S}^{n-1}$ is compact, every distribution $\nu$ on $\mathbb{S}^{n-1}$ is of \emph{finite order}, that is,
there exist $k \in \mathbb{N}$ and $C > 0$ such that
\begin{equation} \label{deforderdist}
|\nu(f)| \leq C\, \|f\|_{C^k}
\end{equation}
for every $f \in C^{\infty}(\mathbb{S}^{n-1})$. The \emph{order} of $\nu$ is the smallest $k$ such that (\ref{deforderdist}) holds.

Examples of distributions on a general smooth manifold $M$ are smooth densities on $M$. Therefore, distributions are often also called generalized densities and
$C^{-\infty}(|\Lambda|(M))$ is used to denote the space of distributions on $M$. However, the choice of a Riemannian metric on $M$ induces an isomorphism between
the space of generalized densities and the space of generalized functions on $M$, usually denoted by $C^{-\infty}(M)$ (cf.\ \textbf{\cite{guilleminsterberg}}). Throughout this article, when $M = \mathbb{S}^{n-1}$,we always make use of this identification and, thus, write $C^{-\infty}(\mathbb{S}^{n-1})$ for the space of distributions on
$\mathbb{S}^{n-1}$ equipped with the topology of weak convergence. The canonical bilinear pairing on $C^{\infty}(\mathbb{S}^{n-1}) \times C^{-\infty}(\mathbb{S}^{n-1})$ will be denoted by
$\langle \,\cdot\,,\,\cdot\, \rangle$.

Since every (signed) Borel measure $\mu$ on $\mathbb{S}^{n-1}$ defines a distribution
$\nu_{\mu}$ by
\[\langle f, \nu_{\mu} \rangle = \int_{\mathbb{S}^{n-1}} f(u)\,d\mu(u), \qquad f \in C^{\infty}(\mathbb{S}^{n-1}),  \]
we will use the continuous linear injection $\mu \mapsto \nu_{\mu}$, to identify $\mathcal{M}(\mathbb{S}^{n-1})$ with a subspace of
$C^{-\infty}(\mathbb{S}^{n-1})$. Clearly, this subspace consists precisely of the distributions of order $0$.
In the same way, the spaces $C^{\infty}(\mathbb{S}^{n-1})$, $C(\mathbb{S}^{n-1})$, and $L^2(\mathbb{S}^{n-1})$ can be
identified with (dense) subspaces of $C^{-\infty}(\mathbb{S}^{n-1})$ and we have
\begin{equation} \label{inclusion1742}
C^{\infty}(\mathbb{S}^{n-1}) \subseteq C(\mathbb{S}^{n-1}) \subseteq L^2(\mathbb{S}^{n-1})
\subseteq \mathcal{M}(\mathbb{S}^{n-1}) \subseteq C^{-\infty}(\mathbb{S}^{n-1}).
\end{equation}

The natural action of $\mathrm{SO}(n)$ on $C^{-\infty}(\mathbb{S}^{n-1})$ is defined as follows: For
$\vartheta \in \mathrm{SO}(n)$ and $\nu \in C^{-\infty}(\mathbb{S}^{n-1})$, we set
\begin{equation} \label{actionson}
\langle f,\vartheta \nu \rangle = \langle \vartheta^{-1} f,\nu \rangle, \qquad f \in C^{\infty}(\mathbb{S}^{n-1}).
\end{equation}
Note that if $\nu$ is a measure on $\mathbb{S}^{n-1}$, then $\vartheta \nu$ is just the image measure of $\nu$ under the rotation $\vartheta$ and
that (\ref{actionson}) is also consistent with the left regular representation of $\mathrm{SO}(n)$ on the spaces $C^{\infty}(\mathbb{S}^{n-1})$, $C(\mathbb{S}^{n-1})$, and $L^2(\mathbb{S}^{n-1})$.
We also remark that the action $\mathrm{SO}(n) \times C^{-\infty}(\mathbb{S}^{n-1}) \rightarrow C^{-\infty}(\mathbb{S}^{n-1})$ is \emph{not} continuous
if $C^{-\infty}(\mathbb{S}^{n-1})$ is endowed with the weak topology, but it is continuous in the strong topology, that is, the topology of uniform convergence on bounded sets.

In order to extend the orthogonal projection $\pi_k: L^2(\mathbb{S}^{n-1}) \rightarrow
\mathcal{H}_k^n$ to $C^{-\infty}(\mathbb{S}^{n-1})$, note that $\pi_k$ is self-adjoint. In particular,
$\langle f,\pi_k g \rangle = \langle \pi_k f, g \rangle$ for all $f \in C^{\infty}(\mathbb{S}^{n-1})$ and $g \in L^2(\mathbb{S}^{n-1})$. In view of (\ref{inclusion1742}), it is therefore consistent
to define $\pi_k \nu$ for $\nu \in C^{-\infty}(\mathbb{S}^{n-1})$ as the distribution given by
\[\langle f, \pi_k \nu \rangle = \langle \pi_k f, \nu \rangle, \qquad f \in C^{\infty}(\mathbb{S}^{n-1}).  \]
From this, it follows (cf.\ \textbf{\cite[\textnormal{p.\ 38}]{morimoto98}}) that indeed $\pi_k \nu \in \mathcal{H}_k^n$ for every $k \in \mathbb{N}$.

\vspace{0.2cm}

Next let us discuss the convolution of functions and measures on $\mathbb{S}^{n-1}$.
Recall that the convolution $\sigma \ast \mu$ of signed measures $\sigma, \mu$ on $\mathrm{SO}(n)$ can be defined by
\[\int_{\mathrm{SO}(n)}\!\!\! f(\vartheta)\, d(\sigma \ast \mu)(\vartheta)=\int_{\mathrm{SO}(n)}\! \int_{\mathrm{SO}(n)}\!\!\! f(\eta \theta)\,d\sigma(\eta)\,d\mu(\theta), \qquad f \in C(\mathrm{SO}(n)).   \]
In other words, $\sigma \ast \mu$ is the pushforward of the product measures $\sigma \otimes \mu$ by the
group multiplication $m: \mathrm{SO}(n) \times \mathrm{SO}(n) \rightarrow \mathrm{SO}(n)$, that is, $\sigma \ast \mu = m_*(\sigma \otimes \mu)$.
Since $\mathrm{SO}(n)$ is compact, this definition can be readily extended to distributions
by replacing the product measure with the tensor product of distributions (see, e.g., \textbf{\cite[\textnormal{p.\ 128}]{hoermander83}}).

The identification of $\mathbb{S}^{n-1}$ with the homogeneous space $\mathrm{SO}(n)/\mathrm{SO}(n-1)$
can now be used to identify $C^{-\infty}(\mathbb{S}^{n-1})$ with right $\mathrm{SO}(n-1)$ invariant distributions on $\mathrm{SO}(n)$. Using this
correspondence, the convolution of distributions on $\mathrm{SO}(n)$ induces a convolution product on $C^{-\infty}(\mathbb{S}^{n-1})$ as follows:
Let $\pi: \mathrm{SO}(n) \rightarrow \mathbb{S}^{n-1}$, $\pi(\eta) = \eta \bar{e}$, be the canonical projection.
Then the convolution of distributions $\delta, \nu \in C^{-\infty}(\mathbb{S}^{n-1})$ is defined by
\[\delta \ast \nu = \pi_* m_*(\pi^*\delta \otimes \pi^*\nu),   \]
where $\pi_*$ and $\pi^*$ denote the pushforward and pullback by $\pi$, respectively.

\pagebreak

The convolution product defined in this way has the following well-known continuity property (see, e.g., \textbf{\cite[\textnormal{Chapter 6}]{hoermander83}}).

\begin{lem} \label{weakcontconv} If $\nu_j \in C^{-\infty}(\mathbb{S}^{n-1})$, $j \in \mathbb{N}$, converge weakly to $\nu \in C^{-\infty}(\mathbb{S}^{n-1})$, then
$\lim_{j \rightarrow \infty} \delta \ast \nu_j = \delta \ast \nu$ weakly for every $\delta \in C^{-\infty}(\mathbb{S}^{n-1})$.
\end{lem}

Zonal measures play an essential role for spherical convolution. For later use we state here explicitly the expressions
for the convolution of a function $h \in C(\mathbb{S}^{n-1})$ and a measure $\sigma \in \mathcal{M}(\mathbb{S}^{n-1})$
with a \emph{zonal} measure $\mu \in \mathcal{M}(\mathbb{S}^{n-1})$ and a \emph{zonal} function $f \in C(\mathbb{S}^{n-1})$, respectively:
\begin{equation} \label{zonalconv}
(h \ast \mu)(\bar{\eta}) = \int_{\mathbb{S}^{n-1}} h(\eta u)\,d\mu(u) \quad \mbox{and} \quad (\sigma \ast f)(\bar{\eta}) = \int_{\mathbb{S}^{n-1}} f(\eta^{-1} u)\,d\sigma(u),
\end{equation}
where for $\eta \in \mathrm{SO}(n)$, we write $\pi(\eta) = \bar{\eta} \in \mathbb{S}^{n-1}$.

From (\ref{zonalconv}) one can deduce several properties of the convolution with zonal measures,
for example, that the convolution of zonal functions and measures is Abelian and that for every $\vartheta \in \mathrm{SO}(n)$, we have $(\vartheta \sigma) \ast \mu = \vartheta(\sigma \ast \mu)$.
Moreover, from the identification of the zonal measure $\mu$ on $\mathbb{S}^{n-1}$ with a measure on $[-1,1]$ and the well-known Funk--Hecke Theorem, it follows (cf.\ \textbf{\cite{Schu06a}})
that the Fourier expansion of $\sigma \ast \mu$ is given by
\begin{equation} \label{funkheckgen}
\sigma \ast \mu \sim \sum_{k=0}^{\infty} a_k^n[\mu]\,\pi_k\sigma,
\end{equation}
where the numbers $a_k^n[\mu]$ are defined by
\begin{equation*}
a_k^n[\mu] = \omega_{n-1} \int_{-1}^1 P_k^n(t)\,(1-t^2)^{\frac{n-3}{2}}\,d\mu(t).
\end{equation*}

Like the convolution of functions on $\mathbb{R}^n$, spherical convolution can be used to
approximate a given function or distribution on $\mathbb{S}^{n-1}$ by smooth functions. To this end, let
$B_j(\bar{e})$, $j \in \mathbb{N}$, denote the open geodesic ball of radius $\frac{1}{j}$ centered at $\bar{e} \in \mathbb{S}^{n-1}$. A
sequence of non-negative \emph{zonal} functions $h_j \in C^{\infty}(\mathbb{S}^{n-1})$, $j \in \mathbb{N}$, is called a {\it spherical approximate
identity} if for each $j$,
\begin{equation} \label{sphapprprop}
\int_{\mathbb{S}^{n-1}}h_j(u)\,du=1 \qquad \mbox{and} \qquad \mathrm{supp}\,h_j \subseteq B_j(\bar{e}).
\end{equation}
For a proof of the following auxiliary result, we refer to \textbf{\cite{grinbergzhang99}} or \textbf{\cite[\textnormal{Chapter 6}]{rudin}}.

\begin{lem} \label{approxid} If $h_j \in C^{\infty}(\mathbb{S}^{n-1})$, $j \in \mathbb{N}$, is a spherical approximate identity, then
\begin{enumerate}
\item[(i)] $\lim_{j \rightarrow \infty} g \ast h_j = g$ uniformly for every $g \in C(\mathbb{S}^{n-1})$;
\item[(ii)] $\lim_{j \rightarrow \infty} \nu \ast h_j = \nu$ weakly for every $\nu \in C^{-\infty}(\mathbb{S}^{n-1})$.
\end{enumerate}
\end{lem}

In the final part of this section, we turn to the Christoffel problem and its solution by Berg \textbf{\cite{cberg}}. First recall that spherical harmonics are eigenfunctions of the Laplacian $\Delta_\mathbb{S}$, more precisely, for $Y_k \in \mathcal{H}_k^n$,
\begin{equation} \label{deltasmult}
\Delta_{\mathbb{S}} Y_k = -k(k + n - 2)\,Y_k.
\end{equation}
Like the orthogonal projection $\pi_k$, the Laplacian $\Delta_{\mathbb{S}}$ is self-adjoint. Thus, it is consistent to define
$\Delta_{\mathbb{S}}\nu$ for $\nu \in C^{-\infty}(\mathbb{S}^{n-1})$ as the distribution given by
\[\langle f,\Delta_{\mathbb{S}} \nu \rangle = \langle \Delta_{\mathbb{S}} f, \nu \rangle, \qquad f \in C^{\infty}(\mathbb{S}^{n-1}).  \]
In particular, by (\ref{inclusion1742}), $\Delta_{\mathbb{S}}$ can now be applied to
support functions of not necessarily smooth convex bodies. This is important for us, since the first-order area measure $S_1(K,\cdot)$ of $K \in \mathcal{K}^n$
and its support function $h_K$ are related by a \emph{linear} differential operator $\Box_n$ in the following way:
\begin{equation} \label{boxhks1}
S_1(K,\cdot) = h_K +  \frac{1}{n-1}\Delta_{\mathbb{S}}h_K =: \Box_n h_K.
\end{equation}

From (\ref{deltasmult}) and the definition of $\Box_n$, it follows that for $f \in C^{\infty}(\mathbb{S}^{n-1})$, the spherical harmonic
expansion of $\Box_n f$ is given by
\begin{equation} \label{boxnmult}
\Box_n f \sim \sum_{k=0}^{\infty} \frac{(1-k)(k+n-1)}{n-1} \pi_kf.
\end{equation}
Hence, the kernel of $\Box_n: C^{\infty}(\mathbb{S}^{n-1}) \rightarrow C^{\infty}(\mathbb{S}^{n-1})$ is given by $\mathcal{H}_1^n$, that is, it consists precisely of the
restrictions of linear functions on $\mathbb{R}^n$ to $\mathbb{S}^{n-1}$. In the following let
\[C_{\mathrm{o}}^{\infty}(\mathbb{S}^{n-1}) = \{f \in C^{\infty}(\mathbb{S}^{n-1}): \pi_1 f = 0\}  \]
and define $C_{\mathrm{o}}^{-\infty}(\mathbb{S}^{n-1})$ in the same way. Then $\Box_n: C_{\mathrm{o}}^{\infty}(\mathbb{S}^{n-1}) \rightarrow C_{\mathrm{o}}^{\infty}(\mathbb{S}^{n-1})$
is an $\mathrm{SO}(n)$ equivariant isomorphism of topological vector spaces.

The classical Christoffel problem consists in finding necessary and sufficient conditions for a Borel measure on $\mathbb{S}^{n-1}$ to be
the first-order area measure of a convex body. A solution was obtained by Berg \textbf{\cite{cberg}} by providing an explicit inversion
formula for the operator $\Box_n$. He proved that for every $n \geq 2$ there exists a uniquely
determined $C^{\infty}$ function $g_n$ on $(-1,1)$ such that the associated zonal function $\breve{g}_n(u)=g_n(u\cdot \bar{e})$ is in $L^1(\mathbb{S}^{n-1})$ and
\begin{equation} \label{multgn}
a_1^n[g_n] = 0, \quad \qquad a_k^n[g_n] = \frac{n-1}{(1-k)(k+n-1)}, \quad k \neq 1.
\end{equation}
It follows from (\ref{boxnmult}), (\ref{funkheckgen}), and (\ref{multgn}) that for every $f \in C_\mathrm{o}^{\infty}(\mathbb{S}^{n-1})$,
\begin{equation} \label{boxninverse}
f = (\Box_n f) \ast \breve{g}_n.
\end{equation}

In the final section, we need a generalization of (\ref{boxninverse}) that can be deduced from \linebreak \textbf{\cite[\textnormal{Theorem 4.3}]{goodeyweil2}} and was independently proved in \textbf{\cite{BPSW2014}}:
For every $j \in \{2, \ldots, n\}$, the integral transform
\[\mathrm{T}_{g_j}: C_{\mathrm{o}}^{\infty}(\mathbb{S}^{n-1}) \rightarrow C_{\mathrm{o}}^{\infty}(\mathbb{S}^{n-1}), \quad f \mapsto f \ast \breve{g}_j,  \]
is an isomorphism. We denote by $\Box_j: C_{\mathrm{o}}^{\infty}(\mathbb{S}^{n-1}) \rightarrow C_{\mathrm{o}}^{\infty}(\mathbb{S}^{n-1})$ the inverse of $\mathrm{T}_{g_j}$.

\vspace{1cm}

\centerline{\large{\bf{ \setcounter{abschnitt}{3}
\arabic{abschnitt}. Smooth and Generalized Valuations}}}

\reseteqn \alpheqn \setcounter{theorem}{0}

\vspace{0.6cm}

We now turn to the background material on translation invariant scalar and convex body valued valuations. In particular, we recall the definitions of
smooth and generalized (spherical) valuations as well as the Alesker--Poincar\'e duality map.

If $G$ is a group of affine transformations on $\mathbb{R}^n$, a
valuation $\phi$ is called \linebreak $G$-invariant if $\phi(gK) =
\phi(K)$ for all $K \in \mathcal{K}^n$ and every $g \in G$. Let $\mathbf{Val}$ denote the vector space of \emph{continuous} translation invariant
scalar valued valuations on $\mathbb{R}^n$. \linebreak It was first proved by McMullen \textbf{\cite{McMullen77}} that
\begin{equation} \label{mcmullen}
\mathbf{Val} = \bigoplus_{0 \leq i \leq n} \mathbf{Val}_i,
\end{equation}
where $\mathbf{Val}_i \subseteq \mathbf{Val}$ denotes the subspace of valuations (homogeneous) of degree $i$.

Recall that a map $\Phi: \mathcal{K}^n \rightarrow \mathcal{K}^n$ is called a \emph{Minkowski valuation} if
\[\Phi K + \Phi L = \Phi(K \cup L) + \Phi(K \cap L)  \]
whenever $K \cup L$ is convex and addition on $\mathcal{K}^n$ is Minkowski addition. We denote by $\mathbf{MVal}$ the set of
continuous translation invariant Minkowski valuations, and we write $\mathbf{MVal}_i$, $0 \leq i \leq n$, for its subset of Minkowski valuations
of degree $i$.

More general than Minkowski valuations, we also consider valuations with values in $C(\mathbb{S}^{n-1})$,
that is, maps $F: \mathcal{K}^n \rightarrow C(\mathbb{S}^{n-1})$, $K \mapsto F_K$, such that
\[F_K + F_L = F_{K \cup L} + F_{K \cap L} \]
whenever $K \cup L$ is convex. Let $\mathbf{CVal}$ denote the vector space of all such valuations which are continuous and translation invariant and, as before,
let $\mathbf{CVal}_i$, $0 \leq i \leq n$, denote its subspace of valuations of degree $i$.

Note that any Minkowski valuation $\Phi \in \mathbf{MVal}$ induces a valuation $F^{\Phi} \in \mathbf{CVal}$ by $F^{\Phi}_K = h_{\Phi K}$,
and that $\Phi$ is $\mathrm{SO}(n)$ equivariant if and only if $F^{\Phi}$ is $\mathrm{SO}(n)$ equivariant. Using the map $\Phi \mapsto F^{\Phi}$,
the set $\mathbf{MVal}$ can be identified with an infinite dimensional convex cone in $\mathbf{CVal}$.

Clearly, a valuation $F \in \mathbf{CVal}$ is uniquely determined by the family of valuations $\varphi_u \in \mathbf{Val}$, $u \in \mathbb{S}^{n-1}$, defined by
$\varphi_u(K) = F_K(u)$ for $K \in \mathcal{K}^n$. If in addition $F$ is $\mathrm{SO}(n)$ equivariant, then for $\eta \in \mathrm{SO}(n)$ and $\bar{\eta} = \eta \bar{e} \in \mathbb{S}^{n-1}$,
\[\varphi_{\bar{\eta}}(K) = F_K(\eta \bar{e}) = (\eta^{-1}F_K)(\bar{e}) = F_{\eta^{-1}K}(\bar{e}) = \varphi_{\bar{e}}(\eta^{-1}K).  \]

\pagebreak

\noindent Hence, an $\mathrm{SO}(n)$ equivariant $F \in \mathbf{CVal}$ is uniquely determined by a single $\mathrm{SO}(n - 1)$ invariant valuation $\varphi_{\bar{e}} \in \mathbf{Val}$.
In fact, there is a one-to-one correspondence between the subspace of $\mathrm{SO}(n - 1)$ invariant valuations in $\mathbf{Val}$ and the subspace
of $\mathrm{SO}(n)$ equivariant valuations in $\mathbf{CVal}$. This observation leads to the following:

\vspace{0.3cm}

\noindent {\bf Definition} \emph{Suppose that $F \in \mathbf{CVal}$ is $\mathrm{SO}(n)$ equivariant. The $\mathrm{SO}(n - 1)$
invariant valuation $\varphi \in \mathbf{Val}$, defined by}
\[\varphi(K) = F_K(\bar{e}), \qquad K \in \mathcal{K}^n,  \]
\emph{is called the \emph{associated real valued valuation} of $F \in \mathbf{CVal}$.}

\vspace{0.3cm}

The following collection of examples and results on \emph{homogeneous} valuations will be useful for later reference.

\begin{expls} \label{exps} \end{expls}

\vspace{-0.3cm}

\begin{enumerate}
\item[(a)] It is a trivial fact that $\mathbf{Val}_0$ is one-dimensional and spanned by the Euler characteristic $V_0$.
(Recall that $V_0(K) = 1$ for every $K \in \mathcal{K}^n$.) Using this observation, it follows that $\Phi_0 \in \mathbf{MVal}_0$
if and only if there exists an $L_0 \in \mathcal{K}^n$ such that $\Phi_0K = L_0$ for every $K \in \mathcal{K}^n$.
If $\Phi_0$ is also $\mathrm{SO}(n)$ equivariant, then $L_0 = c_0 B$, where $B$ denotes the Euclidean unit ball in $\mathbb{R}^n$ and $c_0 \geq 0$.

Hadwiger \textbf{\cite[\textnormal{p.\ 79}]{hadwiger51}} proved that also $\mathbf{Val}_n$ is one-dimensional and spanned by the ordinary volume $V_n$.
From this one can easily deduce that $\Phi_n \in \mathbf{MVal}_n$ if and only if there exists an $L_n \in \mathcal{K}^n$ such that $\Phi_nK = L_nV(K)$ for every $K \in \mathcal{K}^n$.
If $\Phi_n$ is also $\mathrm{SO}(n)$ equivariant, then $L_n = c_n B$ for some $c_n \geq 0$.

\item[(b)] It was first proved by Spiegel \textbf{\cite{spiegel}} that if $\psi \in \mathbf{Val}_1$, then
\[\psi(K + L) = \psi(K) + \psi(L)  \]
for all $K, L \in \mathcal{K}^n$. Using this Minkowski additivity, a description of valuations in $\mathbf{Val}_1$ was obtained by Goodey and Weil \textbf{\cite{goodeyweil0}}
and refined by Kiderlen \textbf{\cite{kiderlen05}}. In order to state their result, recall that any $f \in C^{\infty}(\mathbb{S}^{n-1})$ can be written
as a difference of two support functions $f = h_{K_f} - h_{r_fB}$,
where $r_f \geq 0$ (see, e.g., \linebreak \textbf{\cite[\textnormal{Lemma 1.7.8}]{schneider93}}). Now for $\psi \in \mathbf{Val}_1$, let $\nu_{\psi} \in C^{-\infty}_{\mathrm{o}}(\mathbb{S}^{n-1})$ be given by
\begin{equation} \label{distnupsi}
\langle f,\nu_{\psi} \rangle = \psi(K_f) - \psi(r_f B), \qquad f \in C^{\infty}(\mathbb{S}^{n-1}),
\end{equation}
which is well defined by the Minkowski additivity of $\psi$. Moreover, since $r_f$ depends continuously
on $f$ in the $C^2$ norm, the distribution $\nu_{\psi} \in C^{-\infty}_{\mathrm{o}}(\mathbb{S}^{n-1})$ is of order at most 2.
This allows one to conclude that if $\psi \in \mathbf{Val}_1$, then there exists a uniquely determined $\nu_{\psi} \in C^{-\infty}_{\mathrm{o}}(\mathbb{S}^{n-1})$
of order at most 2, which can be extended to the vector space spanned by support functions, such that
\[\psi(K) = \langle h_K, \nu_\psi \rangle  \]
for every $K \in \mathcal{K}^n$. Also observe that $\psi \mapsto \nu_{\psi}$ is continuous as a map from $\mathbf{Val}_1$ to $C^{-\infty}_{\mathrm{o}}(\mathbb{S}^{n-1})$.

\pagebreak

Using this description of valuations in $\mathbf{Val}_1$, Kiderlen \textbf{\cite{kiderlen05}} proved that if \linebreak $\Phi_1 \in \mathbf{MVal}_1$ is $\mathrm{SO}(n)$ equivariant, then there exists a uniquely determined \emph{zonal} $\nu_1 \in C^{-\infty}_{\mathrm{o}}(\mathbb{S}^{n-1})$ of order at most 2 such that
\begin{equation} \label{descrip1}
h_{\Phi_1K} = h_K \ast \nu_1
\end{equation}
for every $K \in \mathcal{K}^n$. From (\ref{descrip1}) and a straightforward generalization of (\ref{funkheckgen}),
it follows that for every $\Phi_1 \in \mathbf{MVal}_1$, there exists a uniquely determined sequence of real numbers
$a_k^n[\Phi_1]$, $k \in \mathbb{N}$, such that $a_1^n[\Phi_1] = 0$ and
\begin{equation} \label{minkval1mult}
\pi_k h_{\Phi_1K} =  a_k^n[\Phi_1] \pi_k h_K
\end{equation}
for every $K \in \mathcal{K}^n$ and $k \in \mathbb{N}$. In fact, relation (\ref{minkval1mult}) was already proved by Schneider \textbf{\cite{schneider74}}
in 1974, where he also showed that for every $k \neq 1$,
\begin{equation} \label{schneidest}
 |a_k^n[\Phi_1]| \leq a_0^n[\Phi_1].
\end{equation}
We note that a precise description of the cone of zonal distributions in $C^{-\infty}_{\mathrm{o}}(\mathbb{S}^{n-1})$ which generate a Minkowski
valuation by (\ref{descrip1}) is still open. However, Kiderlen \textbf{\cite{kiderlen05}} showed that this cone contains
all non-negative zonal measures on $\mathbb{S}^{n-1}$. More precisely, if $\mu_1 \in \mathcal{M}_{\mathrm{o}}(\mathbb{S}^{n-1})$ is
zonal and non-negative, then
\[h_{\Psi_1K} = h_K \ast \mu_1, \qquad K \in \mathcal{K}^n,  \]
defines an $\mathrm{SO}(n)$ equivariant Minkowski valuation in $\mathbf{MVal}_1$.

\item[(c)] A classification of continuous translation invariant scalar valued valuations of degree $n - 1$ was obtained
by McMullen \textbf{\cite{McMullen80}}. It states that $\phi \in \mathbf{Val}_{n-1}$ if and only if
there exists a uniquely determined $f \in C_{\mathrm{o}}(\mathbb{S}^{n-1})$ such that
\[\phi(K) = \int_{\mathbb{S}^{n-1}} f(u)\,dS_{n-1}(K,u)  \]
for every $K \in \mathcal{K}^n$.

Applying McMullen's result to associated real valued valuations, the first author \textbf{\cite{Schu06a}} obtained the following
description of $\mathrm{SO}(n)$ equivariant Minkowski valuations of degree $n - 1$:
If $\Phi_{n-1} \in \mathbf{MVal}_{n-1}$ is $\mathrm{SO}(n)$ equivariant, then there exists a uniquely determined \emph{zonal}
$f_{n-1} \in C_{\mathrm{o}}(\mathbb{S}^{n-1})$ such that
\begin{equation} \label{descripn-1}
h_{\Phi_{n-1}K} = S_{n-1}(K,\cdot) \ast f_{n-1}
\end{equation}
for every $K \in \mathcal{K}^n$. As in the case of Minkowski valuations of degree $1$, a precise description of the cone of zonal functions in $C_{\mathrm{o}}(\mathbb{S}^{n-1})$ which
generate a Minkowski valuation by (\ref{descripn-1}) is still open (see \textbf{\cite{Schu06a}} for more information).

\item[(d)] Several important Minkowski valuations arise from data about sections or projections of convex bodies
and are therefore objects of intensive research in geometric tomography (see, e.g., \textbf{\cite{gardner2ed, haberl11, kiderlen05, ludwig02, Ludwig:Minkowski}}).
Of particular interest for us are the {\it normalized mean section operators} $\mathrm{M}_j \in \mathbf{MVal}_{n+1-j}$, $2 \leq j \leq n$, introduced by Goodey and Weil \textbf{\cite{goodeyweil1}}.
In \textbf{\cite[\textnormal{Theorem 4.4}]{goodeyweil2}}, they showed that for $K \in \mathcal{K}^n$,
\begin{equation} \label{genfctmj}
h_{\mathrm{M}_jK} = q_{n,j}\,S_{n+1-j}(K,\cdot) \ast \breve{g}_j,
\end{equation}
where $g_j$ is the $j$th Berg function and the constant $q_{n,j}$ is given by
\[q_{n,j} = \frac{j-1}{2\pi(n+1-j)}\,\frac{\kappa_{j-1}\kappa_{j-2}\kappa_{n-j}}{\kappa_{j-3}\kappa_{n-2}}.\]
Here, $\kappa_i$ is the $i$-dimensional volume of the $i$-dimensional Euclidean unit ball.
\end{enumerate}

\vspace{0.1cm}

A simple consequence of (\ref{mcmullen}) is that the space $\mathbf{Val}$ becomes a Banach space when endowed with
the norm
\[\|\phi \|=\sup \{|\phi(K)|: K \subseteq B\}.  \]
The natural \emph{continuous} action of the general linear group $\mathrm{GL}(n)$ on the Banach space $\mathbf{Val}$ is for $A \in \mathrm{GL}(n)$ given by
\[(A\phi)(K)=\phi(A^{-1}K), \qquad \phi \in \mathbf{Val},\ K \in \mathcal{K}^n.  \]

The notion of smooth vectors of a continuous representation now gives rise to the notion of smooth valuations, first introduced
by Alesker \textbf{\cite{Alesker03}}.

\vspace{0.3cm}

\noindent {\bf Definition} \emph{A valuation $\phi \in \mathbf{Val}$ is called \emph{smooth} if the map $\mathrm{GL}(n) \rightarrow \mathbf{Val}$, defined by $A \mapsto A\phi$, is
infinitely differentiable.}

\vspace{0.3cm}

Note that smooth valuations are precisely the smooth vectors of the natural representation of $\mathrm{GL}(n)$ on $\mathbf{Val}$. We
therefore write $\mathbf{Val}^{\infty}$ for the Fr\'echet space of smooth translation invariant valuations endowed
with the G{\aa}rding topology (see Section 2). We denote the subspace of smooth valuations of degree $i$ by $\mathbf{Val}_i^{\infty}$.
By the general properties of smooth vectors discussed in Section~2, the spaces $\mathbf{Val}_i^{\infty}$ are dense
$\mathrm{GL}(n)$ invariant subspaces of $\mathbf{Val}_i$ and from (\ref{mcmullen}) one can deduce that
\begin{equation*}
\mathbf{Val}^{\infty} = \bigoplus \limits_{0 \leq i \leq n} \mathbf{Val}_i^{\infty}.
\end{equation*}

Since in this article we are mainly concerned with $\mathrm{SO}(n)$ equivariant valuations, the following
result will be useful.

\begin{prop} \label{glnsonsmooth} A valuation $\phi \in \mathbf{Val}$ is smooth if and only if the restricted map $\mathrm{SO}(n) \rightarrow \mathbf{Val}$, $\vartheta \mapsto \vartheta\phi$,
is smooth. Moreover, the G{\aa}rding topologies on $\mathbf{Val}^{\infty}$ induced by the natural representations of $\mathrm{GL}(n)$ and $\mathrm{SO}(n)$ coincide.
\end{prop}

The first statement of Proposition \ref{glnsonsmooth} follows for example from a recent result of Bernstein and Kr\"otz \textbf{\cite[\textnormal{Corollary 3.10}]{bersteinkroetz2014}} and the fact that
the natural representation of $\mathrm{GL}(n)$ on $\mathbf{Val}$ is admissible and of finite length (see \textbf{\cite{Alesker01}}). The second statement is then a simple consequence of
the open mapping theorem.

\vspace{0.1cm}

The one-to-one correspondence between $\mathrm{SO}(n)$ equivariant valuations in $\mathbf{CVal}$ and $\mathrm{SO}(n - 1)$ invariant valuations in $\mathbf{Val}$ now
motivates the following.

\vspace{0.3cm}

\noindent {\bf Definition} \emph{Let $F \in \mathbf{CVal}$ be $\mathrm{SO}(n)$ equivariant and let $\varphi \in \mathbf{Val}$ be the associated real valued valuation of $F$.
We define the \emph{norm} of $F$ by
\begin{equation} \label{normcval}
\|F\| = \sup \{|\varphi(K)|: K \subseteq B\}.
\end{equation}
Furthermore, we call $F$ \emph{smooth} if $\varphi$ is smooth.}

\vspace{0.3cm}

While it is easy to see that McMullen's decomposition result (\ref{mcmullen}) implies that
\begin{equation} \label{cvaldecomp}
\mathbf{CVal} = \bigoplus_{0 \leq i \leq n} \mathbf{CVal}_i,
\end{equation}
it was recently proved by Parapatits and the second author \textbf{\cite{parapwann}} that, in general, a Minkowski valuation $\Phi \in \mathbf{MVal}$
need not be a sum of homogeneous Minkowski valuations $\Phi_i \in \mathbf{MVal}_i$.
However, from (\ref{cvaldecomp}) one can still deduce a
decomposition result for translation invariant Minkowski valuations (cf.\ \textbf{\cite{schnschu}}), which we
state here under the additional assumption of $\mathrm{SO}(n)$ equivariance.

\begin{lem} \label{lemdec} If $\Phi \in \mathbf{MVal}$ is $\mathrm{SO}(n)$ equivariant, then there
exist uniquely determined $c_0, c_n \geq 0$ and $\mathrm{SO}(n)$ equivariant valuations $F_i \in \mathbf{CVal}_i$, $1 \leq i \leq n - 1$ such that
\[h_{\Phi K} = c_0 + \sum \limits_{i=1}^{n-1} F_{i,K} + c_nV(K),  \]
for every $K \in \mathcal{K}^n$. Moreover, if $\Phi$ is smooth, then each $F_i$ is also smooth.
\end{lem}

The vector space of $\mathrm{SO}(n)$ equivariant valuations in $\mathbf{CVal}$ endowed with the norm (\ref{normcval}) becomes a Banach space
in which smooth valuations form a dense subspace. However, it is a priori not clear that an $\mathrm{SO}(n)$ equivariant Minkowski valuation in $\mathbf{MVal}$ can be approximated by smooth ones. We will prove this in Section 5.

\vspace{0.2cm}

It is well known that for any $\phi \in \mathbf{Val}$ and $K \in \mathcal{K}^n$, McMullen's decomposition (\ref{mcmullen}) implies that the function
$t \mapsto \phi(K + tB)$ is a polynomial of degree at most $n$. This, in turn, gives rise to a derivation operator $\Lambda: \mathbf{Val} \rightarrow \mathbf{Val}$, defined by
\begin{equation} \label{derivation1}
(\Lambda \phi)(K) = \left . \frac{d}{dt} \right |_{t=0} \phi(K + tB).
\end{equation}
Using (\ref{derivation1}) it is not difficult to show that if $\phi \in \mathbf{Val}_i$, then $\Lambda \phi \in \mathbf{Val}_{i-1}$, that $\Lambda$ is \linebreak continuous, $\mathrm{SO}(n)$ equivariant, and
that $\Lambda$ maps smooth valuations to smooth ones.

\pagebreak

The advantage of working with smooth translation invariant valuations instead of merely continuous ones is
that the space $\mathbf{Val}^{\infty}$ admits more algebraic structure. For example, the following Hard Lefschetz type theorem for the operator $\Lambda$ was proved
by Alesker \textbf{\cite{Alesker03}} for even and by Bernig and Br\"ocker \textbf{\cite{Bernig07b}} for general valuations.

\begin{theorem}[\!\!\cite{Alesker03, Bernig07b}] \label{hardlef1}
Suppose that $\frac{n}{2} < i \leq n$. Then $\Lambda^{2i-n}: \mathbf{Val}^{\infty}_{i} \rightarrow \mathbf{Val}^{\infty}_{n-i}$ is an isomorphism.
\end{theorem}

Recently, Parapatits and the first author \textbf{\cite{parapschu}} proved that for any $\Phi \in \mathbf{MVal}$, there exist
$\Phi^{(j)} \in \mathbf{MVal}$, where $0 \leq j \leq n$, such that
\[\Phi(K + tB) = \sum_{j=0}^n t^{n-j}\Phi^{(j)}(K)  \]
for every $K \in \mathcal{K}^n$ and $t \geq 0$.
This Steiner type formula shows that the operator $\Lambda$ from (\ref{derivation1}) has a natural analogue for Minkowski valuations $\Lambda:~\mathbf{MVal}~\rightarrow~\mathbf{MVal}$.

\vspace{0.3cm}

\noindent {\bf Definition} \emph{For $\Phi \in \mathbf{MVal}$, define $\Lambda \Phi \in \mathbf{MVal}$ by}
\[ h_{(\Lambda \Phi)(K)}(u) = \left . \frac{d}{dt} \right |_{t=0} h_{\Phi(K + tB)}(u), \qquad u \in \mathbb{S}^{n-1}.  \]

\vspace{0.2cm}

Note that if $\Phi \in \mathbf{MVal}_i$ is $\mathrm{SO}(n)$ equivariant, then
so is $\Lambda\Phi \in \mathbf{MVal}_{i-1}$. Moreover, if $\varphi \in \mathbf{Val}_i$ is the associated real valued valuation of $\Phi$,
then $\Lambda \varphi \in \mathbf{Val}_{i-1}$ is associated with $\Lambda\Phi$. In particular, if $\Phi$ is smooth, then so is $\Lambda \Phi$.

\vspace{0.2cm}

Another important structural property of smooth valuations is the existence of a
\emph{continuous} bilinear product, discovered by Alesker \textbf{\cite{Alesker04a}},
\[\cdot:\mathbf{Val}^{\infty} \times \mathbf{Val}^{\infty} \rightarrow \mathbf{Val}^{\infty}, \quad (\phi,\psi) \mapsto \phi \cdot \psi.\]
Endowed with the Alesker product, $\mathbf{Val}^{\infty}$ becomes an associative and commutative
algebra with unit given by the Euler characteristic which is graded by the degree of homogeneity, that is,
\begin{equation} \label{grading}
\mathbf{Val}_i^{\infty} \cdot \mathbf{Val}_j^{\infty} \subseteq \mathbf{Val}_{i+j}^{\infty}.
\end{equation}

Recall that $\mathbf{Val}_n$ is $1$-dimensional and spanned by $V_n$. If $V_n^* \in \mathbf{Val}_n^*$ is the unique element such that $\langle V_n,V_n^* \rangle = 1$,
then it follows from (\ref{grading}) that for every $0 \leq i \leq n$,
\begin{equation} \label{pairing17}
\langle\, \cdot\,,\cdot \rangle :\mathbf{Val}^{\infty}_i \times \mathbf{Val}^{\infty}_{n-i} \rightarrow \mathbb{R}, \qquad (\phi,\psi) \mapsto \langle \phi \cdot \psi, V_n^* \rangle,
\end{equation}
defines a continuous bilinear pairing between smooth valuations of complementary degree. Moreover, Alesker \textbf{\cite{Alesker04a}} proved that this pairing is \emph{non-degenerate}.

\vspace{0.3cm}

\noindent {\bf Definition} \emph{The space of translation invariant \emph{generalized valuations} of degree \linebreak $i \in \{0, \ldots, n\}$ is defined as the topological dual
\[\mathbf{Val}^{-\infty}_i = (\mathbf{Val}_{n-i}^{\infty})^*\]
endowed with the weak topology.}

\pagebreak

Since the pairing (\ref{pairing17}) is non-degenerate, the \emph{Poincar\'e duality map}, defined by
\begin{equation}
\mathrm{pd}: \mathbf{Val}_i^{\infty} \rightarrow \mathbf{Val}_i^{-\infty}, \qquad \langle \mathrm{pd}\, \phi,\psi \rangle = \langle \phi,\psi \rangle,
\end{equation}
is continuous, injective and has dense image with respect to the weak topology. Moreover,
it follows from \textbf{\cite[\textnormal{Proposition 8.1.2}]{alesker10}} that $\mathrm{pd}$ admits a unique continuous
extension to $\mathbf{Val}_i$. Thus, we will use the Poincar\'e duality map to \emph{identify} the spaces $\mathbf{Val}_i^{\infty}$
and $\mathbf{Val}_i$ with dense subspaces of $\mathbf{Val}_i^{-\infty}$. Hence, we have the inclusions
\[ \mathbf{Val}_i^{\infty} \subseteq \mathbf{Val}_i \subseteq \mathbf{Val}_i^{-\infty}. \]
We also note that the natural action of $\mathrm{SO}(n)$ on $\mathbf{Val}_i^{-\infty}$ is \emph{not} continuous in the weak topology,
but it is continuous when $\mathbf{Val}_i^{-\infty}$ is given the strong topology.

\vspace{0.1cm}

Next we recall a recent result of Alesker, Bernig and the
first author \textbf{\cite{ABS2011}} on the decomposition of the vector spaces of translation invariant (generalized) valuations into $\mathrm{SO}(n)$ irreducible subspaces.

\begin{theorem}[\!\!\cite{ABS2011}] \label{thm_decomposition}
For $0 \leq i \leq n$, the spaces $\mathbf{Val}_i^{\infty}$, $\mathbf{Val}_i$, and $\mathbf{Val}_i^{-\infty}$ are
multiplicity free under the action of $\mathrm{SO}(n)$. Moreover, the highest weights of the $\mathrm{SO}(n)$ irreducible subspaces in either
of them are precisely given by the tuples $(\lambda_1,\ldots,\lambda_{\lfloor n/2 \rfloor})$ satisfying (\ref{heiwei}) and the following additional conditions:
\[(i)\ \lambda_j = 0 \mbox{ for } j > \min\{i,n-i\}; \quad \, (ii)\ |\lambda_j| \neq 1 \mbox{ for } 1 \leq j \leq \lfloor \mbox{$\frac{n}{2}$} \rfloor; \quad\, (iii)\  |\lambda_2| \leq 2.  \]
\end{theorem}

\vspace{0.1cm}

We now use the notion of spherical representations (see Section 2) to
define spherical valuations.

\vspace{0.3cm}

\noindent {\bf Definition} \emph{For $0 \leq i \leq n$, the subspaces $\mathbf{Val}_i^{\mathrm{sph}}$, $\mathbf{Val}_i^{\infty,\mathrm{sph}}$, and $\mathbf{Val}_i^{-\infty,\mathrm{sph}}$
of translation invariant continuous, smooth, and generalized \emph{spherical valuations} of degree $i$ are defined as the closures (w.r.t.\ the respective topologies) of the direct sum of all
$\mathrm{SO}(n)$ irreducible subspaces in $\mathbf{Val}_i$, $\mathbf{Val}_i^{\infty}$, and $\mathbf{Val}_i^{-\infty}$, respectively,
which are spherical with respect to $\mathrm{SO}(n - 1)$.}

\vspace{0.3cm}

Note that, by Theorems \ref{thm_decomposition} and \ref{thmspher} (see also Example \ref{exp2}), $\mathbf{Val}_i^{-\infty,\mathrm{sph}}$ is the \emph{annihilator}
of the closure of the direct sum of all $\mathrm{SO}(n)$ irreducible subspaces in $\mathbf{Val}_{n-i}^{\infty}$ with highest weights not  of the form  $(k,0,\ldots,0)$, $k \in \mathbb{N}$.
As a consequence, every $\mathrm{SO}(n - 1)$ invariant (generalized) valuation in $\mathbf{Val}_i$, $\mathbf{Val}_i^{\infty}$, or $\mathbf{Val}_i^{-\infty}$ is spherical.

\begin{expls} \label{exps2} \end{expls}

\vspace{-0.2cm}

\begin{enumerate}
\item[(a)] It follows from Theorem \ref{thm_decomposition} or Example \ref{exps} (b) that
\[\quad \mathbf{Val}_1 = \mathbf{Val}_1^{\mathrm{sph}}, \qquad \! \quad \mathbf{Val}_1^{\infty} = \mathbf{Val}_1^{\infty,\mathrm{sph}}, \qquad\,  \mathbf{Val}_1^{-\infty} = \mathbf{Val}_1^{-\infty,\mathrm{sph}}.   \]
\item[(b)] It follows from Theorem \ref{thm_decomposition} or Example \ref{exps} (c) that
\[\mathbf{Val}_{n-1} = \mathbf{Val}_{n-1}^{\mathrm{sph}}, \qquad  \mathbf{Val}_{n-1}^{\infty} = \mathbf{Val}_{n-1}^{\infty,\mathrm{sph}}, \qquad  \mathbf{Val}_{n-1}^{-\infty} = \mathbf{Val}_{n-1}^{-\infty,\mathrm{sph}}.   \]
\end{enumerate}

\pagebreak

\centerline{\large{\bf{ \setcounter{abschnitt}{4}
\arabic{abschnitt}. Auxiliary Results about Smooth Valuations}}}

\reseteqn \alpheqn \setcounter{theorem}{0}

\vspace{0.6cm}

In this section we begin with the proof of Theorem \ref{main2} (a).
As a corollary, we obtain a version of Theorem \ref{main1} for \emph{smooth} Minkowski valuations.
We also determine an explicit expression for the pairing (\ref{pairing17}) when one of the valuations is spherical.
This will be needed in the next section to complete the proof of Theorem \ref{main2} (a).

\begin{theorem} \label{smoothsphiso} For $1 \leq i \leq n - 1$, the map $\mathrm{E}_i: C_{\mathrm{o}}^{\infty}(\mathbb{S}^{n-1}) \rightarrow \mathbf{Val}_i^{\infty,\mathrm{sph}}$, defined by
\begin{equation} \label{defei}
(\mathrm{E}_if)(K) = \int_{\mathbb{S}^{n-1}} f(u)\,dS_i(K,u),
\end{equation}
is an $\mathrm{SO}(n)$ equivariant isomorphism of topological vector spaces.
\end{theorem}
{\it Proof.} Clearly, the maps $\mathrm{E}_i: C_{\mathrm{o}}(\mathbb{S}^{n-1}) \rightarrow \mathbf{Val}_i$, given by (\ref{defei}), are
linear and $\mathrm{SO}(n)$ equivariant for every $i \in \{1, \ldots, n - 1\}$. Using the monotonicity of the mixed volumes
$V(K[i];B[n-i])$ with $i$ copies of $K$ and $n-i$ copies of the Euclidean unit ball $B$, we obtain for $K \subseteq B$,
\[|\mathrm{E}_if(K)| \leq \int_{\mathbb{S}^{n-1}} |f(u)| \,dS_i(K,u) \leq nV(K[i];B[n-i])\|f\| \leq n\kappa_n \|f\|, \]
which proves $\|\mathrm{E}_if\| \leq n\kappa_n \|f\|$ and, hence, the continuity of $\mathrm{E}_i$. Thus, by
Lemma~\ref{smoothlem17} and Proposition \ref{glnsonsmooth}, the restrictions $\mathrm{E}_i: C_{\mathrm{o}}^{\infty}(\mathbb{S}^{n-1}) \rightarrow \mathbf{Val}_i^{\infty}$ are well defined and continuous.

Since differences of area measures of order~$i$ of convex bodies in $\mathcal{K}^n$ are dense in the
set of all signed finite Borel measures on $\mathbb{S}^{n-1}$ with centroid at the origin (see, e.g., \textbf{\cite[\textnormal{p.\ 477}]{schneider93}})),
the maps $\mathrm{E}_i$ are also injective. Consequently, by Schur's lemma and Example \ref{exp2}, $\mathrm{E}_i(\mathcal{H}_k^n)$, $k \neq 1$, is an
$\mathrm{SO}(n)$ irreducible subspace of $\mathbf{Val}_i^{\infty}$ of highest weight $(k,0,\ldots,0)$. By the definition of the spaces $\mathbf{Val}_i^{\infty,\mathrm{sph}}$,
it follows that $\mathbf{U}_i := \mathrm{E}_i(C_{\mathrm{o}}^{\infty}(\mathbb{S}^{n-1}))$ is a dense subspace of $\mathbf{Val}_i^{\infty,\mathrm{sph}}$. By the open mapping theorem, it remains to show that $\mathbf{U}_i$ is closed.

First, let $i = n - 1$. The result of McMullen \textbf{\cite{McMullen80}} discussed in Example \ref{exps} (c) implies that
the map $\mathrm{E}_{n-1}: C_{\mathrm{o}}(\mathbb{S}^{n-1}) \rightarrow \mathbf{Val}_i^{\mathrm{sph}}$
is a continuous bijection and hence, by the open mapping theorem, an isomorphism of Banach spaces. The assertion
for $i = n - 1$ is now an immediate consequence of Lemma~\ref{smoothlem17} and Proposition \ref{glnsonsmooth}.

Next, recall that the area measures satisfy the Steiner type formula
\[S_i(K + tB,\cdot) = \sum_{j=0}^i t^{i-j} {i \choose j} S_j(K,\cdot)  \]
for every $K \in \mathcal{K}^n$ and $t \geq 0$. Thus, for $f \in C_{\mathrm{o}}^{\infty}(\mathbb{S}^{n-1})$ and $i \geq 2$, we have
\begin{equation} \label{ablsph}
(\Lambda \mathrm{E}_if)(K) =  \left . \frac{d}{dt} \right |_{t=0}  \int_{\mathbb{S}^{n-1}} f(u)\,dS_i(K + tB,u) = i(\mathrm{E}_{i-1}f)(K).
\end{equation}
In particular, the restriction of the derivation operator $\Lambda$ to $\mathbf{U}_i$ is injective for $2 \leq i \leq n - 1$ and $\Lambda(\mathbf{U}_i) = \mathbf{U}_{i-1}$.
Consequently, the restriction of the linear $\mathrm{SO}(n)$ equivariant map $\Lambda$ to  $\mathrm{cl}\,\mathbf{U}_i = \mathbf{Val}_i^{\infty,\mathrm{sph}}$ is injective as well.
By Theorem \ref{hardlef1}, $\Lambda^{n-2}: \mathbf{Val}_{n-1}^{\infty} \rightarrow \mathbf{Val}_1^{\infty}$ is an $\mathrm{SO}(n)$ equivariant isomorphism.
But $\mathbf{U}_{n-1} = \mathbf{Val}_{n-1}^{\infty}$ by what we have proved above. Therefore,
\[\mathbf{U}_1 = \Lambda^{n-2}(\mathbf{U}_{n-1}) = \mathbf{Val}_1^{\infty}.  \]
In particular, $\mathbf{U}_1 = \Lambda^{i-1}(\mathbf{U}_i)$ is closed.
Hence, $\Lambda^{i-1}(\mathbf{U}_i) \subseteq \Lambda^{i-1}(\mathrm{cl}\,\mathbf{U}_i) \subseteq \Lambda^{i-1}(\mathbf{U}_i)$ that is,
\[
\Lambda^{i-1}(\mathbf{U}_i) = \Lambda^{i-1}(\mathrm{cl}\,\mathbf{U}_i).
\]
Since $\Lambda^{i-1}$ is injective on $\mathrm{cl}\,\mathbf{U}_i = \mathbf{Val}_i^{\infty,\mathrm{sph}}$, we conclude that $\mathbf{U}_i$ is closed and, therefore,
$\mathbf{U}_i = \mathbf{Val}_i^{\infty,\mathrm{sph}}$. \hfill $\blacksquare$

\vspace{0.3cm}

Using Lemma \ref{lemdec} and Theorem \ref{smoothsphiso}, we can now prove the following Hadwiger type result for \emph{smooth} Minkowski valuations.

\begin{koro} \label{mainsmooth} If $\Phi: \mathcal{K}^n \rightarrow \mathcal{K}^n$ is a smooth Minkowski valuation which is translation invariant and $\mathrm{SO}(n)$ equivariant,
then there exist uniquely determined $c_0, c_n \geq 0$ and $\mathrm{SO}(n-1)$ invariant $f_i \in C_{\mathrm{o}}^{\infty}(\mathbb{S}^{n-1})$, $1 \leq i \leq n - 1$,  such that
\begin{equation}
h_{\Phi K} = c_0 + \sum_{i=1}^{n-1} S_i(K,\cdot) \ast f_i + c_n V_n(K)
\end{equation}
for every $K \in \mathcal{K}^n$.
\end{koro}
{\it Proof.} By Lemma \ref{lemdec}, it is sufficient to show that for every $\mathrm{SO}(n)$ equivariant smooth $F_i \in \mathbf{CVal}_i$, $1 \leq i \leq n - 1$,
there exists a uniquely determined $\mathrm{SO}(n-1)$ invariant $f_i \in C_{\mathrm{o}}^{\infty}(\mathbb{S}^{n-1})$, such that for every $K \in \mathcal{K}^n$,
\begin{equation} \label{toshowsmooth}
F_{i,K} = S_i(K,\cdot) \ast f_i.
\end{equation}
To this end, let $\varphi_i \in \mathbf{Val}_i^\infty$ denote the associated real valued valuation of $F_i$ and recall that, by definition, $\varphi_i$ is $\mathrm{SO}(n - 1)$ invariant.
It follows that $\varphi_i \in \mathbf{Val}_i^{\infty,\mathrm{sph}}$. Thus, by Theorem \ref{smoothsphiso}, there exists
a uniquely determined $f_i \in C_{\mathrm{o}}^{\infty}(\mathbb{S}^{n-1})$ such that
\[\varphi_i(K) = \int_{\mathbb{S}^{n-1}} f_i(u)\,dS_i(K,u).   \]
Moreover, the $\mathrm{SO}(n - 1)$ invariance of $\varphi_i$ implies that also $f_i$ is $\mathrm{SO}(n - 1)$ invariant. Hence, by the definition of $\varphi_i$ and (\ref{zonalconv}), we obtain
\[F_{i,K}(\bar{\eta}) = \varphi_i(\eta^{-1}K) =  \int_{\mathbb{S}^{n-1}} f_i(\eta^{-1}u)\,dS_i(K,u) = (S_i(K,\cdot) \ast f_i)(\bar{\eta}), \]
where for $\eta \in \mathrm{SO}(n)$, we set, as before, $\bar{\eta} = \eta \bar{e} \in \mathbb{S}^{n-1}$. \hfill $\blacksquare$

\vspace{0.3cm}

Corollary \ref{mainsmooth} under the additional assumption that the Minkowski valuation $\Phi$ is \emph{even}
was recently obtained by the authors \textbf{\cite{SchuWan13}} using a different approach.

\pagebreak

Finally, we require a generalization of a formula of Bernig and Hug \textbf{\cite{bernighug2015}} for the pairing (\ref{pairing17}) of spherical valuations.
To this end, let $a: \mathbb{S}^{n-1} \rightarrow \mathbb{S}^{n-1}$ denote the antipodal map, given by, $a(u) = -u$, $u \in \mathbb{S}^{n-1}$,
and recall from Example \ref{exps} (b) that any $\psi \in \mathbf{Val}_1$ determines a unique $\nu_{\psi} \in C^{-\infty}_{\mathrm{o}}(\mathbb{S}^{n-1})$, defined by (\ref{distnupsi}).

\begin{prop} \label{poincsph} Let $1 \leq i \leq n - 1$. For $\phi_i \in \mathbf{Val}_i$ and $f \in C^{\infty}_{\mathrm{o}}(\mathbb{S}^{n-1})$, we have
\begin{equation} \label{yoda17}
 \langle \phi_i,\mathrm{E}_{n-i}f \rangle = \frac{(n-i)!}{(n-1)!}\, \langle f \circ a, \nu_{\Lambda^{i-1}\phi_i} \rangle.
\end{equation}
\end{prop}
{\it Proof.} Since both pairings in (\ref{yoda17}) are jointly continuous and bilinear, we may assume that $f \in \mathcal{H}_k^n$ for some $k \in \mathbb{N}$, $k \neq 1$, and
that $\phi_i$ belongs to an $\mathrm{SO}(n)$ irreducible subspace $\Gamma_{\lambda} \subseteq \mathbf{Val}_i$ of highest weight $\lambda = (\lambda_1,\ldots,\lambda_{\lfloor n/2\rfloor})$.
In particular, $\phi_i$ is smooth.

Next, note that $\Lambda^{i-1} \phi_i \in \mathbf{Val}_1^{\infty} = \mathbf{Val}_1^{\infty,\mathrm{sph}}$ (cf.\ Examples \ref{exps2}). Therefore, by Theorem \ref{smoothsphiso}, there exists a smooth function $h \in C^{\infty}_{\mathrm{o}}(\mathbb{S}^{n-1})$ (in fact, $h \in \mathcal{H}_m^n$ for some $m \in \mathbb{N}$) such that
\begin{equation} \label{luke17}
(\Lambda^{i-1}\phi_i)(K) = i! \int_{\mathbb{S}^{n-1}} h(u)\,dS_1(K,u).
\end{equation}
The normalizing coefficient $i!$ is chosen for convenience as will become clear below.

Since the pairing (\ref{pairing17}) is biinvariant under the (simultaneous) action of $\mathrm{SO}(n)$ and the spaces $\mathcal{H}_k^n$ are
self-dual as $\mathrm{SO}(n)$ modules, the restriction of the Poincar\'e duality map to $\Gamma_{\lambda}$ defines a linear
$\mathrm{SO}(n)$ equivariant map from $\Gamma_{\lambda}$ to $\mathcal{H}_k^n$, that is,
\[ \mathrm{pd}|_{\Gamma_{\lambda}} \in \mathrm{Hom}_{\mathrm{SO}(n)}(\Gamma_{\lambda},\mathcal{H}_k^n). \]
Since both $\Gamma_{\lambda}$ and $\mathcal{H}_k^n$ are $\mathrm{SO}(n)$ irreducible, it follows from Schur's lemma
that $\mathrm{pd}|_{\Gamma_{\lambda}}$ and, thus, the left hand side of (\ref{yoda17}) can only be non-zero when $\Gamma_{\lambda}$ and $\mathcal{H}_k^n$ are
isomorphic, that is, when $(\lambda_1,\ldots,\lambda_{\lfloor n/2\rfloor}) = (k,0,\ldots,0)$. Similarly, the right hand side of (\ref{yoda17})
can only be non-zero if $(\lambda_1,\ldots,\lambda_{\lfloor n/2\rfloor}) = (k,0,\ldots,0)$. We may therefore assume that $\phi_i$ is spherical.
But if $\phi_i \in \mathbf{Val}_i^{\infty,\mathrm{sph}}$, then, by (\ref{luke17}) and (\ref{ablsph}), we have
\[\phi_i(K) = \int_{\mathbb{S}^{n-1}} h(u)\,dS_i(K,u) = \mathrm{E}_ih.   \]
In this case, it follows from \textbf{\cite[\textnormal{Proposition 4.11}]{bernighug2015}} that
\begin{equation} \label{hansolo17}
\langle \mathrm{E}_ih,\mathrm{E}_{n-i}f \rangle = \frac{(n-i)!i!}{(n-1)!} \int_{\mathbb{S}^{n-1}} h(u)\, \Box_nf(-u)\,du.
\end{equation}
Finally, definition (\ref{distnupsi}), (\ref{luke17}), and (\ref{boxhks1}) yield
\[\langle \phi_i,\mathrm{E}_{n-i}f \rangle = \frac{(n-i)!}{(n-1)!}\, \langle f \circ a, \nu_{\Lambda^{i-1}\phi_i} \rangle.\]
\hfill $\blacksquare$

\pagebreak

\centerline{\large{\bf{ \setcounter{abschnitt}{5}
\arabic{abschnitt}. Proof of the Main Results}}}

\reseteqn \alpheqn \setcounter{theorem}{0}

\vspace{0.6cm}

We are now in a position to complete the proofs of Theorems \ref{main1} and \ref{main2}.
We also discuss a more precise version of Theorem \ref{main1} for homogeneous Minkowski valuations in dimensions $n \leq 4$.
At the end of the section we include an approximation result for continuous Minkowski valuations by smooth ones.

We begin with the following slightly more precise version of Theorem \ref{main2} (a).

\begin{theorem} \label{eistern} For $1 \leq i \leq n - 1$, the isomorphism $\mathrm{E}_i: C_{\mathrm{o}}^{\infty}(\mathbb{S}^{n-1}) \rightarrow \mathbf{Val}_i^{\infty,\mathrm{sph}}$
admits a unique extension by continuity in the weak topologies to an isomorphism
\[\widetilde{\mathrm{E}}_i: C_{\mathrm{o}}^{-\infty}(\mathbb{S}^{n-1}) \rightarrow \mathbf{Val}_i^{-\infty,\mathrm{sph}}.   \]
Moreover, the diagram
\[
\begin{xy}
 \xymatrix{
   C^{\infty}_{\mathrm{o}}(\mathbb{S}^{n-1}) \, \ar[rr]^{\textnormal{\normalsize $\mathrm{E}_i$}} \ar@{->}[dd] \ar@_{(->}[dd] &     & \, \mathbf{Val}_i^{\infty,\,\mathrm{sph}} \ar@{->}[dd] \ar@_{(->}[dd]_{\textnormal{\normalsize $\mathrm{pd}$}}  \\
       &     &     \\
      C^{-\infty}_{\mathrm{o}}(\mathbb{S}^{n-1}) \,    \ar[rr]^{\textnormal{\normalsize $\widetilde{\mathrm{E}}_i$}} &     & \, \mathbf{Val}_i^{-\infty,\,\mathrm{sph}}
  }
\end{xy}
\]
commutes and the vertical maps have dense image.
\end{theorem}
{\it Proof.} First recall that $\mathbf{Val}_i^{-\infty,\mathrm{sph}}$ is the annihilator of the subspace spanned by all
$\mathrm{SO}(n)$ irreducible subspaces of $\mathbf{Val}_{n-i}^{\infty}$ which are non-spherical. Hence, using Theorem \ref{smoothsphiso},
we can define a map $\widetilde{\mathrm{E}}_i: C_{\mathrm{o}}^{-\infty}(\mathbb{S}^{n-1}) \rightarrow \mathbf{Val}_i^{-\infty,\mathrm{sph}}$ by
\begin{equation} \label{defeistern}
\langle \widetilde{\mathrm{E}}_i\nu, \phi_{n-i} \rangle = \frac{i!}{(n-1)!}\, \langle \nu, (\Box_n \circ a^* \circ \mathrm{E}_{1}^{-1}\circ \Lambda^{n-i})\phi_{n-i} \rangle,
\end{equation}
where $\phi_{n-i} \in \mathbf{Val}_{n-i}^{\infty}$ and $a^*$ denotes the pullback by the antipodal map. From (\ref{hansolo17}), it follows that $\widetilde{\mathrm{E}}_i$ continuously extends $\mathrm{E}_i$.

Since the differential operator $\Box_n:  C^{\infty}_{\mathrm{o}}(\mathbb{S}^{n-1}) \rightarrow  C^{\infty}_{\mathrm{o}}(\mathbb{S}^{n-1})$ is an isomorphism, it follows from
Theorem~\ref{smoothsphiso} that $\widetilde{\mathrm{E}}_i\nu = 0$ implies $\nu = 0$, that is, $\widetilde{\mathrm{E}}_i$ is injective. In order to prove
that $\widetilde{\mathrm{E}}_i$ is surjective, let $\xi \in \mathbf{Val}_i^{-\infty,\mathrm{sph}}$ be given and note that
$\xi \circ \mathrm{E}_{n-i} \in C^{-\infty}_{\mathrm{o}}(\mathbb{S}^{n-1})$. If we put
\[\nu = \frac{(n-1)!}{(n-i)!i!}\,\xi \circ \mathrm{E}_{n-i} \circ a^* \circ \Box_n^{-1} \in C_{\mathrm{o}}^{-\infty}(\mathbb{S}^{n-1}),   \]
then, by (\ref{defeistern}),
\[\langle \widetilde{\mathrm{E}}_i\nu,\mathrm{E}_{n-i}f \rangle = \langle \xi,\mathrm{E}_{n-i}f \rangle  \]
for all $f \in C_{\mathrm{o}}^{\infty}(\mathbb{S}^{n-1})$, that is, $\widetilde{\mathrm{E}}_i \nu = \xi$. Clearly, the map $\widetilde{\mathrm{E}}_i^{-1}$ thus defined is continuous. \hfill $\blacksquare$

\pagebreak

In the next lemma, which is crucial for the proof of Theorem \ref{main2} (b), and in all that follows, the letter $C$ will denote a constant
that can be different from one line to the next and that depends only on the dimension $n$.

\begin{lem} \label{maintech} There exists a constant $C > 0$ such that
\[\|f\|_{C^2} \leq C\,\|\Box_n f\|_{C^0}.   \]
for every $\mathrm{SO}(n - 1)$ invariant $f \in C^2_{\mathrm{o}}(\mathbb{S}^{n-1})$.
\end{lem}
{\it Proof.} For arbitrary but fixed $q \in \mathbb{R}$, we consider the linear differential operator $\mathrm{D}_q: C^2(\mathbb{S}^{n-1}) \rightarrow C(\mathbb{S}^{n-1})$, defined by
\[\mathrm{D}_q f = \Delta_{\mathbb{S}} f + qf. \]
Note that $\mathrm{D}_q$ is $\mathrm{SO}(n)$ equivariant and that $\mathrm{D}_{n-1} = (n-1)\,\Box_n$. Moreover, by (\ref{deltasmult}), the operator $\mathrm{D}_q$
is injective for every $q \neq k(k + n - 2)$, $k \in \mathbb{N}$. If $q = k(k + n - 2)$ for some $k \in \mathbb{N}$, then the kernel of $\mathrm{D}_q$ is given by $\mathcal{H}_k^n$.

First, we show that there exists a constant $C > 0$ such that for $k = 1, 2$,
\begin{equation} \label{impinequ}
|\nabla^k f|_0 := \max\limits_{\mathbb{S}^{n-1}}|\nabla^k f| \leq C(\|f\|_{C^0} + \|\mathrm{D}_qf\|_{C^0})
\end{equation}
for every $\mathrm{SO}(n - 1)$ invariant $f \in C^2(\mathbb{S}^{n-1})$ and, therefore, by (\ref{defcknorm}),
\begin{equation} \label{first1}
\|f\|_{C^2} \leq C(\|f\|_{C^0} + \|\mathrm{D}_qf\|_{C^0}).
\end{equation}
Using the cylindrical coordinates (\ref{cylindcoord}) and expressions (\ref{normnabf}) and (\ref{normnab2f}),
we see that in order to prove (\ref{impinequ}), it suffices to prove that
\begin{equation} \label{toprove1}
|\partial_t f|_0 := \sup\limits_{(-1,1)} \left | \frac{\partial f}{\partial t} \right | \leq C(\|f\|_{C^0} + \|\mathrm{D}_qf\|_{C^0})
\end{equation}
and
\begin{equation} \label{toprove2}
|(1-t^2)\partial_{tt}f|_0 := \sup\limits_{(-1,1)} \left |(1-t^2)\frac{\partial^2 f}{\partial t^2} \right | \leq C(\|f\|_{C^0} + \|\mathrm{D}_qf\|_{C^0})
\end{equation}
for every $\mathrm{SO}(n - 1)$ invariant $f \in C^2(\mathbb{S}^{n-1})$. But since $(1-t^2)\partial_{tt}f = \Delta_{\mathbb{S}}f + (n-1)t\partial_tf$ by (\ref{lapcylininv}), it follows that (\ref{toprove2}) is actually an immediate consequence of (\ref{toprove1}) and the definition of $\mathrm{D}_q$. Thus, in order to prove (\ref{impinequ}), we only have to show that (\ref{toprove1}) holds
for every $\mathrm{SO}(n - 1)$ invariant $f \in C^2(\mathbb{S}^{n-1})$.

Let $f \in C^2(\mathbb{S}^{n-1})$ now be an arbitrary but fixed $\mathrm{SO}(n - 1)$ invariant function. Since $|\nabla f|$ and $|\nabla^2 f|$ are bounded on $\mathbb{S}^{n-1}$, it follows from (\ref{normnabf}) and (\ref{normnab2f}) that $\partial_t f$ is bounded on $(-1,1)$. Assume that $|\partial_t f|$ attains its maximum at $t_0 \in (-1,1)$. Since, by (\ref{lapcylininv}),
\[\Delta_{\mathbb{S}}f = (1-t^2)\partial_{tt}f - (n-1)t\partial_tf = (1-t^2)^{1-(n-1)/2}\partial_t\left ((1-t^2)^{(n-1)/2}\partial_t f \right ), \]
it follows from the definition of $\mathrm{D}_q$ that
\begin{eqnarray*}
(1-t_0^2)^{(n-1)/2}\partial_tf(t_0) &  = & \int_{-1}^{t_0}\partial_t \left ( (1-t^2)^{(n-1)/2}\partial_t f(t) \right )\, dt \\
& = & \int_{-1}^{t_0} (1-t^2)^{(n-1)/2-1}(\mathrm{D}_qf(t)-qf(t))\,dt
\end{eqnarray*}
and, hence,
\begin{equation} \label{batman}
|(1-t_0^2)^{(n-1)/2}\partial_tf(t_0)| \leq C(\|f\|_{C^0} + \|\mathrm{D}_qf\|_{C^0}).
\end{equation}
This shows that we may assume that $|t_0| \geq \alpha$ for some fixed $\alpha > 0$, otherwise (\ref{batman}) implies (\ref{toprove1}).
But since $\partial_{tt} f(t_0) = 0$, we conclude from (\ref{lapcylininv}) that
\[-(n-1)t_0\partial_tf(t_0) = \mathrm{D}_qf(t_0) -qf(t_0)  \]
which also yields $|\partial_t f(t_0)| \leq C(\|f\|_{C^0} + \|\mathrm{D}_qf\|_{C^0})$. Hence, we have shown that
\begin{equation} \label{zwischen17}
|\partial_t f|_0 \leq \max \left \{C(\|f\|_{C^0} + \|\mathrm{D}_qf\|_{C^0}), \limsup \limits_{t \rightarrow \pm 1} |\partial_t f(t)| \right \}
\end{equation}
and it remains to bound $\limsup_{t \rightarrow \pm 1} |\partial_t f(t)|$ in terms of $\|f\|_{C^0}$ and $\|\mathrm{D}_qf\|_{C^0}$. In order to do this, note that, by (\ref{lapcylininv}),
$\partial_t f$ is a \emph{bounded} solution on $(-1,1)$ of the differential equation
\[y'(t) - \frac{(n-1)t}{1-t^2}y(t) = \frac{\mathrm{D}_qf(t)-qf(t)}{1-t^2}.   \]
All solutions of this equation are given by
\[y(t) = (1-t^2)^{-(n-1)/2} \left (\int_{-1}^t \frac{\mathrm{D}_qf(s)-qf(s)}{(1-s^2)^{1-(n-1)/2}}\,ds + c \right ),  \]
where $c \in \mathbb{R}$. Since $\partial_t f$ is bounded, we must have
\begin{equation} \label{ironman}
\partial_tf(t) = (1-t^2)^{-(n-1)/2} \int_{-1}^t \frac{\mathrm{D}_qf(s)-qf(s)}{(1-s^2)^{1-(n-1)/2}}\,ds
\end{equation}
and
\begin{equation} \label{hulk}
\int_{-1}^1 \frac{\mathrm{D}_qf(s)-qf(s)}{(1-s^2)^{1-(n-1)/2}}\,ds = 0.
\end{equation}
Consequently,
\begin{eqnarray*}
\limsup \limits_{t \rightarrow 1} |\partial_t f(t)| & = & \limsup \limits_{t \rightarrow 1} (1-t^2)^{-(n-1)/2} \left | \int_{t}^1 \frac{\mathrm{D}_qf(s)-qf(s)}{(1-s^2)^{1-(n-1)/2}}\,ds \right | \\
 & \leq & \limsup \limits_{t \rightarrow 1} \frac{1-t}{1-t^2}  \|\mathrm{D}_qf-qf\|_{C^0} \leq C(\|f\|_{C^0} + \|\mathrm{D}_qf\|_{C^0}).
\end{eqnarray*}
Similarly, we obtain $\limsup_{t \rightarrow -1} |\partial_t f(t)| \leq  C(\|f\|_{C^0} + \|\mathrm{D}_qf\|_{C^0})$ which, by (\ref{zwischen17}),
completes the proof of (\ref{toprove1}) and thus of (\ref{impinequ}) and (\ref{first1}).

\pagebreak

Next, assume that $q < 0$. Then the maximum principle implies that there exists a $C > 0$ such that
\[\|f\|_{C^0} \leq C\|\mathrm{D}_q f\|_{C^0}  \]
for every $f \in C^2(\mathbb{S}^{n-1})$. Consequently, we obtain from (\ref{first1}) that
\[\|f\|_{C^2} \leq C\|\mathrm{D}_q f\|_{C^0}.   \]
for every $\mathrm{SO}(n-1)$ invariant $f \in C^2(\mathbb{S}^{n-1})$. From this, it follows that $\mathrm{D}_q$ is injective, has dense image and that
\[\mathrm{D}_q^{-1}: C(\mathbb{S}^{n-1})^{\mathrm{SO}(n-1)} \rightarrow C^2(\mathbb{S}^{n-1})^{\mathrm{SO}(n-1)} \hookrightarrow C(\mathbb{S}^{n-1})^{\mathrm{SO}(n-1)}  \]
exists and is bounded. Here, $C(\mathbb{S}^{n-1})^{\mathrm{SO}(n-1)}$ denotes the Banach subspace of all $\mathrm{SO}(n-1)$ invariant functions in $C(\mathbb{S}^{n-1})$ and
$C^2(\mathbb{S}^{n-1})^{\mathrm{SO}(n-1)}$ is defined similarly. Moreover, the Arzel\`a--Ascoli theorem implies that $\mathrm{D}_q^{-1}$ is \emph{compact}.

Now, choose an $m > n - 1$ and put $q = n - m - 1 < 0$. Applying the Fredholm alternative (see, e.g., \textbf{\cite[\textnormal{Theorem 5.3}]{gilbardtrudinger}})
to the compact operator $\mathrm{D}_q^{-1}: C_{\mathrm{o}}(\mathbb{S}^{n-1})^{\mathrm{SO}(n-1)} \rightarrow C_{\mathrm{o}}(\mathbb{S}^{n-1})^{\mathrm{SO}(n-1)}$ yields that either
\begin{equation} \label{fastfertig1}
f + m\mathrm{D}_q^{-1}f = 0
\end{equation}
has a non-trivial solution $f \in C_{\mathrm{o}}(\mathbb{S}^{n-1})^{\mathrm{SO}(n-1)}$ or
\begin{equation} \label{fastfertig2}
f + m\mathrm{D}_q^{-1}f = \mathrm{D}_q^{-1}h
\end{equation}
has a solution for every $h \in C_{\mathrm{o}}(\mathbb{S}^{n-1})^{\mathrm{SO}(n-1)}$. In the latter case, the operator $(\mathrm{Id}+m\mathrm{D}_q^{-1})^{-1}$ is bounded.
However, since
\[\mathrm{D}_q(f + m\mathrm{D}_q^{-1}f)= \Delta_{\mathbb{S}}f + (n-1)f = 0   \]
for $f \in C_{\mathrm{o}}(\mathbb{S}^{n-1})^{\mathrm{SO}(n-1)}$ implies that $f = 0$,
equation (\ref{fastfertig1}) has \emph{no} non-trivial solution in $C_{\mathrm{o}}(\mathbb{S}^{n-1})^{\mathrm{SO}(n-1)}$ and thus (\ref{fastfertig2})
is solvable for every $h \in C_{\mathrm{o}}(\mathbb{S}^{n-1})^{\mathrm{SO}(n-1)}$, that is, $h = \Delta_{\mathbb{S}}f + (n-1)f$
is solvable for every $h \in C_{\mathrm{o}}(\mathbb{S}^{n-1})^{\mathrm{SO}(n-1)}$ and
\[\|f\|_{C^0} = \|(\mathrm{Id}+m\mathrm{D}_q^{-1})^{-1}\mathrm{D}_q^{-1}h\|_{C^0} \leq C \|\mathrm{D}_q^{-1}h\|_{C^0} \leq C \|h\|_{C^0} = C\|\mathrm{D}_{n-1}f\|_{C^0}.  \]
Combining this with (\ref{first1}) for the case $q = n - 1$ and recalling that $\mathrm{D}_{n-1} = (n - 1)\,\Box_n$, completes the proof of the lemma. \hfill $\blacksquare$

\vspace{0.3cm}

We remark, that Lemma \ref{maintech} without the assumption of $\mathrm{SO}(n-1)$ invariance does not hold in general.

\vspace{0.1cm}

Using Lemma \ref{maintech} and Proposition \ref{poincsph}, we can now complete the proof of Theorem \ref{main2} (b).

\vspace{0.3cm}

\noindent {\it Proof of Theorem \ref{main2} (b).} Let $\phi_i \in \mathbf{Val}_i^{\mathrm{SO}(n - 1)}$, $1 \leq i \leq n - 1$, and recall
that every $\mathrm{SO}(n - 1)$ invariant valuation is spherical. Hence, using the Poincar\'e duality map, we can identify
$\phi_i$ with a generalized valuation from $\mathbf{Val}_i^{-\infty,\mathrm{sph}}$. By Proposition \ref{poincsph},
\[\langle \phi_i,\mathrm{E}_{n-i}f \rangle = \frac{(n-i)!}{(n-1)!}\, \langle f \circ a, \nu_{\Lambda^{i-1}\phi_i} \rangle  \]
for $f \in C^{\infty}_{\mathrm{o}}(\mathbb{S}^{n-1})$.
Since $\Lambda^{i-1}\phi_i$ is $1$-homogeneous, $\nu_{\Lambda^{i-1}\phi_i} \in C^{-\infty}_{\mathrm{o}}(\mathbb{S}^{n-1})$ is of order at most $2$ (cf.\ Example \ref{exps} (b)). Hence, $\phi_i \circ \mathrm{E}_{n-i} \in C^{-\infty}_{\mathrm{o}}(\mathbb{S}^{n-1})$ defines an $\mathrm{SO}(n-1)$ invariant distribution of order at most $2$.

At the same time, by Theorem \ref{eistern}, $\phi_i = \widetilde{\mathrm{E}}_i\gamma$ for some uniquely determined $\gamma \in C^{-\infty}_{\mathrm{o}}(\mathbb{S}^{n-1})$ and
since $\phi_i$ is $\mathrm{SO}(n-1)$ invariant, so is $\gamma$. We want to show that $\gamma$ is of order $0$ and, thus, in fact a measure.
To this end, first note that, by (\ref{defeistern}),
\[\phi_i \circ \mathrm{E}_{n-i} = \frac{(n-i)!i!}{(n-1)!}\, \gamma \circ \Box_n \circ a^*.  \]
Now, for $f \in C^{\infty}_{\mathrm{o}}(\mathbb{S}^{n-1})$, let
\[\overline{f} = \int_{\mathrm{SO}(n-1)} \vartheta f\,d\vartheta = \delta_{\bar{e}} \ast f \]
denote the $\mathrm{SO}(n-1)$-rotational symmetral of $f$. Clearly, we have $\|\overline{f}\,\|_{C^0} \leq \|f\|_{C^0}$.
Moreover, it is not difficult to show (cf.\ \textbf{\cite[\textnormal{Theorem 6.30}]{rudin}}) that the $\mathrm{SO}(n-1)$ invariance of $\gamma$ implies $\gamma(f) = \gamma(\overline{f})$.
Consequently, using Lemma \ref{maintech}, we obtain
\[|\gamma(f)| = |\gamma(\overline{f})| = C |(\phi_i \circ \mathrm{E}_{n-i} \circ a^*)(\Box_n^{-1}\overline{f})| \leq C\|\Box_n^{-1}\overline{f}\|_{C^2} \leq C\|\overline{f}\,\|_{C^0} \leq C\|f\|_{C^0},  \]
that is, $\gamma$ is of order $0$ and therefore a measure.

In the case $i = n - 1$, it follows from the result of McMullen \textbf{\cite{McMullen80}}, described in Example \ref{exps} (c), that, in fact,
$\phi_i \in \widetilde{\mathrm{E}}_i(C_{\mathrm{o}}(\mathbb{S}^{n-1}))$. \hfill $\blacksquare$

\vspace{0.3cm}

In the same way Theorem \ref{smoothsphiso} implies Corollary \ref{mainsmooth}, we can use Theorem \ref{main2} (b) and an approximation argument to deduce
Theorem \ref{main1}.

\vspace{0.3cm}

\noindent {\it Proof of Theorem \ref{main1}.} By Lemma \ref{lemdec}, we have to show that for every $\mathrm{SO}(n)$ equivariant $F_i \in \mathbf{CVal}_i$, $1 \leq i \leq n - 1$,
there exist uniquely determined $\mathrm{SO}(n-1)$ invariant measures $\mu_i \in \mathcal{M}_{\mathrm{o}}(\mathbb{S}^{n-1})$, $1 \leq i \leq n - 2$, and
an $\mathrm{SO}(n - 1)$ invariant function $f_{n-1} \in C_{\mathrm{o}}(\mathbb{S}^{n-1})$, such that for $1 \leq i \leq n - 2$,
\begin{equation} \label{case1nmin2}
F_{i,K} = S_i(K,\cdot) \ast \mu_i
\end{equation}
and
\begin{equation} \label{casenmin1}
F_{n-1,K} = S_{n-1}(K,\cdot) \ast f_{n-1}
\end{equation}
for every $K \in \mathcal{K}^n$.

\pagebreak

Since (\ref{casenmin1}) can be proved, using Theorem \ref{main2} (b), in exactly the same way that (\ref{toshowsmooth}) was
deduced from Theorem \ref{smoothsphiso}, we only explain the proof of (\ref{case1nmin2}) here. \linebreak
First, let $K \in \mathcal{K}^n$ be such that $h_K \in C^{\infty}(\mathbb{S}^{n-1})$ and that $K$ has positive curvature.
Then the area measure $S_i(K,\cdot)$ of $K$ is absolutely continuous with respect to spherical Lebesgue measure with a smooth density function
$s_i(K,\cdot) \in C^{\infty}_{\mathrm{o}}(\mathbb{S}^{n-1})$ (see, e.g., \textbf{\cite[\textnormal{Chapter 2.5}]{schneider93}}).
We want to show that if $\varphi_i \in \mathbf{Val}_i^{\mathrm{sph}}$ denotes the $\mathrm{SO}(n - 1)$ invariant associated real valued valuation of $F_i$, $1 \leq i \leq n - 2$,
then there exists a uniquely determined $\mathrm{SO}(n - 1)$ invariant $\mu_i \in \mathcal{M}_{\mathrm{o}}(\mathbb{S}^{n-1})$ such that
\begin{equation} \label{prf1717}
\varphi_i(K) = \int_{\mathbb{S}^{n-1}} s_i(K,u)\,d\mu_i(u).
\end{equation}
To this end, note that, by Theorem \ref{main2}, there exists a uniquely determined $\mathrm{SO}(n - 1)$ invariant $\mu_i \in \mathcal{M}_{\mathrm{o}}(\mathbb{S}^{n-1})$
such that $\varphi_i = \widetilde{\mathrm{E}}_i\mu_i$. Moreover, it follows from a result of Bernig and Faifman \textbf{\cite[\textnormal{p.\ 11}]{bernigfaifman15}}
that
\begin{equation} \label{bernigfaif1}
\varphi_i(K) = \langle \varphi_i,\psi_{n-i}^K \rangle,
\end{equation}
where $\psi_{n-i}^K \in \mathbf{Val}_{n-i}^{\infty}$ is given by the mixed volume
\[\psi_{n-i}^K(L) = {n \choose i} V(L[n-i],-K[i]).  \]
Now, let $f_{i,j} \in C^{\infty}_{\mathrm{o}}(\mathbb{S}^{n-1})$, $j \in \mathbb{N}$, be a sequence of smooth functions which converges weakly to $\mu_i$.
Then, by (\ref{bernigfaif1}) and Proposition \ref{poincsph}, we have
\[\varphi_i(K) = \lim_{j \rightarrow \infty} \langle \widetilde{\mathrm{E}}_i f_{i,j},\psi_{n-i}^K \rangle = \frac{i!}{(n-1)!} \lim_{j \rightarrow \infty} \langle f_{i,j} \circ a,\nu_{\Lambda^{n-i-1}\psi_{n-i}^K} \rangle.  \]
Using the definitions of $\psi_{n-i}^K$ and $\Lambda$ it is not difficult to show that
\[(\Lambda^{n-i-1}\psi_{n-i}^K)(L) = \frac{n!}{i!} V(L,B[n-i-1],-K[i]) = \frac{(n-1)!}{i!}\int_{\mathbb{S}^{n-1}} h_L(u)dS_i(-K,u).   \]
Thus, using $dS_i(-K,u) = s_i(K,-u)\,du$ and the definition of $\nu_{\Lambda^{n-i-1}\psi_{n-i}^K}$, we obtain
\[\varphi_i(K) = \lim_{j \rightarrow \infty} \int_{\mathbb{S}^{n-1}} s_i(K,u)f_{i,j}(u)\,du = \int_{\mathbb{S}^{n-1}}  s_i(K,u)\,d\mu_i(u) \]
which completes the proof of (\ref{prf1717}).

From the definition of $\varphi_i$, (\ref{prf1717}), and (\ref{zonalconv}), we now obtain
\[F_{i,K}(\bar{\eta}) = \varphi_i(\eta^{-1}K) =  \int_{\mathbb{S}^{n-1}} s_i(K,\eta u)\,d\mu_i(u) = (s_i(K,\cdot) \ast \mu_i)(\bar{\eta}).\]
Since both sides of this equation depend continuously on $K$, (\ref{case1nmin2}) follows from the fact
that convex bodies with smooth support functions and positive curvature are dense in $\mathcal{K}^n$. \hfill $\blacksquare$

\pagebreak

The following consequence of Theorem \ref{main1} for homogeneous Minkowski valuations
includes a slight improvement for dimensions $n \leq 4$ which we deduce from the existence of the derivation
operator $\Lambda: \mathbf{MVal} \rightarrow \mathbf{MVal}$ and the estimate (\ref{schneidest}).

\begin{koro} \label{homogcoro} Let $\Phi_i: \mathcal{K}^n \rightarrow \mathcal{K}^n$ be a continuous, translation invariant, and
$\mathrm{SO}(n)$ equivariant Minkowski valuation of degree $i \in \{0, \ldots, n\}$.
\begin{enumerate}
\item[(i)] If $i = 0$, then $\Phi_0K = c_0B$ for some $c_0 \geq 0$ and every $K \in \mathcal{K}^n$.
\item[(ii)] If $1 \leq i \leq n - 2$, then there exists a uniquely determined $\mathrm{SO}(n - 1)$ invariant $\mu_i \in \mathcal{M}_{\mathrm{o}}(\mathbb{S}^{n-1})$
such that $h_{\Phi_iK} = S_i(K,\cdot) \ast \mu_i$ for every $K \in \mathcal{K}^n$.
\item[(iii)] If $i=n - 1$, then there exists a uniquely determined $\mathrm{SO}(n - 1)$ invariant $f_{n-1} \in C_{\mathrm{o}}(\mathbb{S}^{n-1})$
such that $h_{\Phi_{n-1}K} = S_{n-1}(K,\cdot) \ast f_{n-1}$ for every $K \in \mathcal{K}^n$.
\item[(iv)] If $i = n$, then $\Phi_nK = c_nV_n(K)B$ for some $c_n \geq 0$ and every $K \in \mathcal{K}^n$.
\end{enumerate}
Moreover, if $n = 3$ or $n = 4$, then the measures $\mu_i$, $i = 1, 2$, from (ii) are absolutely continuous with respect to
spherical Lebesgue measure with densities in $L^2_{\mathrm{o}}(\mathbb{S}^{n-1})$.
\end{koro}
{\it Proof.} The statements (i)--(iv) are direct consequences of Theorem \ref{main1}, so we only have to prove
the absolute continuity of the measures $\mu_i$, $i = 1, 2$, for $n \leq 4$. To this end, first note that
$\Lambda^{i-1}\Phi_i \in \mathbf{MVal}_1$ is $\mathrm{SO}(n)$ equivariant and that if $\varphi_i \in \mathbf{Val}_i$
is the associated real valued valuation of $\Phi_i$, then $\Lambda^{i-1}\varphi_i \in \mathbf{Val}_1$ is associated with $\Lambda^{i-1}\Phi_i$.
Thus, it follows easily from (ii), (\ref{ablsph}), (\ref{boxhks1}) and the fact that multiplier transformations commute that for every $K \in \mathcal{K}^n$,
\begin{equation} \label{prfhomimprove}
h_{\Lambda^{i-1}\Phi_iK} = i!\, S_1(K,\cdot) \ast \mu_i = i!\,h_K \ast \Box_n\mu_i.
\end{equation}
Hence, the distribution determined by $\Lambda^{i-1}\varphi_i \in \mathbf{Val}_1$ (cf.\ Example \ref{exps} (b)) is given by $i!\,\Box_n\mu_i$. Since $\mu_i$
is $\mathrm{SO}(n-1)$ invariant, so is $i!\,\Box_n\mu_i$ and it follows from (\ref{expzonal}) and (\ref{boxnmult}) that the Fourier expansion of $i!\,\Box_n\mu_i$ is given by
\[i!\,\Box_n\mu_i \sim i!\sum_{k=0}^{\infty} \frac{N(n,k)}{\omega_n}\,\frac{(1-k)(k+n-1)}{n-1}\,a_k^n[\mu_i]\,P_k^n(\,.\cdot\bar{e}).  \]
Therefore, using (\ref{prfhomimprove}), (\ref{funkheckgen}), and (\ref{schneidest}), it follows that there exists an absolute constant $C > 0$ such that for every $k \geq 2$,
\begin{equation} \label{estimate17}
|a_k^n[\mu_i]| \leq C\,\frac{i!(n-1)}{(k-1)(k+n-1)}.
\end{equation}
But, since $\left( \frac{N(n,k)}{\omega_n}\right )^{1/2}P_k^n(\,.\cdot\bar{e})$ forms an orthonormal sequence in $L^2(\mathbb{S}^{n-1})$ (see, e.g., \textbf{\cite[\textnormal{p.\ 84}]{groemer96}}) and, by (\ref{nnk}), $N(n,k) = \mathrm{O}(k^{n-2})$ as $k \rightarrow \infty$, we see that
\[\mu_i \sim \sum_{k=0}^{\infty} \frac{N(n,k)}{\omega_n}\,a_k^n[\mu_i]\,P_k^n(\,.\cdot\bar{e}).  \]
converges in $L^2(\mathbb{S}^{n-1})$ as long as $n \leq 4$. \hfill $\blacksquare$

\pagebreak

Corollary \ref{homogcoro} (iii) was previously obtained by the first author \textbf{\cite{Schu06a}} as already explained in Example \ref{exps} (c).
The case $i = 1$ of Corollary \ref{homogcoro} (ii) can be reformulated as follows (cf.\ (\ref{prfhomimprove})): There exists a uniquely determined
$\mathrm{SO}(n - 1)$ invariant $\mu_1 \in \mathcal{M}_{\mathrm{o}}(\mathbb{S}^{n-1})$ such that for every $K \in \mathcal{K}^n$,
\begin{equation} \label{kidimprov}
h_{\Phi_1K} = h_K \ast \Box_n\mu_1.
\end{equation}
Comparing (\ref{kidimprov}) with the corresponding result (\ref{descrip1}) of Kiderlen \textbf{\cite{kiderlen05}}, shows that we have slightly improved (\ref{descrip1}) by proving that
the $\mathrm{SO}(n-1)$ invariant distribution $\nu_1 \in C_{\mathrm{o}}^{-\infty}(\mathbb{S}^{n-1})$ of order at most 2 determined by $\Phi_1$ is always of the form $\nu_1 = \Box_n\mu_1$ for some
$\mathrm{SO}(n - 1)$ invariant $\mu_1 \in \mathcal{M}_{\mathrm{o}}(\mathbb{S}^{n-1})$.

Note that the estimate (\ref{estimate17}) is not strong enough to deduce that $\mu_i$ is absolutely continuous in higher dimensions,
as can be seen for example from the spherical Radon (or Minkowski-Funk) transform,
$\mathrm{R}: C(\mathbb{S}^{n-1}) \rightarrow C(\mathbb{S}^{n-1})$, defined by
\[\mathrm{R}f = f \ast \mu_{\mathbb{S}^{n-2}}, \]
where $\mu_{\mathbb{S}^{n-2}} \in \mathcal{M}(\mathbb{S}^{n-1})$ is uniformly concentrated on $\mathbb{S}^{n-1} \cap \bar{e}^{\bot}$. Clearly, $\mu_{\mathbb{S}^{n-2}}$ is $\mathrm{SO}(n-1)$ invariant but not absolutely continuous with respect to spherical Lebesgue measure. However, $|a_k^n[\mu_{\mathbb{S}^{n-2}}]| = \mathrm{O}(k^{1-n/2})$ as $k \rightarrow \infty$ (see, \textbf{\cite[\textnormal{Lemma 3.4.8}]{groemer96}}).

Finally we remark that Corollary \ref{homogcoro} (ii) does not leave much room for improvement since
the zonal functions $\breve{g}_n$ associated with Berg's functions are not continuous on $\mathbb{S}^{n-1}$ for $n \geq 3$ and
they do not lie in $L^2_{\mathrm{o}}(\mathbb{S}^{n-1})$ but merely in $L^1_{\mathrm{o}}(\mathbb{S}^{n-1})$ for $n \geq 5$. However,
they are generating functions of the normalized mean section operators $\mathrm{M}_j$ as described in Example \ref{exps} (d).

\vspace{0.2cm}

We conclude this section with an approximation result of continuous Minkowski valuations by smooth ones which generalizes a
corresponding result for \emph{even} Minkowski valuations of the first author \textbf{\cite{Schu09}} and will be useful in the last section.

\begin{koro} Every continuous translation invariant and $\mathrm{SO}(n)$
equivariant Minkowski valuation can be approximated
uniformly on compact subsets of $\mathcal{K}^n$ by smooth
translation invariant and $\mathrm{SO}(n)$ equivariant
Minkowski valuations.
\end{koro}
{\it Proof.} Let $\Phi \in \mathbf{MVal}$ be $\mathrm{SO}(n)$ equivariant and let
\begin{equation}
h_{\Phi K} = c_0 + \sum_{i=1}^{n-2} S_i(K,\cdot) \ast \mu_i + S_{n-1}(K,\cdot) \ast f_{n-1} + c_n V_n(K)
\end{equation}
be the convolution representation of $\Phi$ according to Theorem \ref{main1}.
We define a sequence $\Phi^j \in \mathbf{MVal}$, $j \in \mathbb{N}$, of $\mathrm{SO}(n)$ equivariant Minkowski valuations by
\[h_{\Phi^jK} = h_{\Phi K} \ast h_j, \qquad K \in \mathcal{K}^n,   \]
where $h_j$, $m \in \mathbb{N}$, is a spherical approximate identity. Note that since $h_j \geq 0$, $\Phi^j$ is well defined
by the result of Kiderlen \textbf{\cite{kiderlen05}} described at the end of Example \ref{exps} (b). Using (\ref{sphapprprop}) and the $\mathrm{SO}(n)$ equivariance of $\Phi$,
it is easy to show that $\Phi^j$ converges to $\Phi$ on compact subsets (cf. the proof of \textbf{\cite[\textnormal{Theorem 6.5}]{Schu09}}).

\pagebreak

It remains to show that the Minkowski valuations $\Phi^j$ are smooth, that is, the associated real valued valuations $\varphi^j \in \mathbf{Val}^{\mathrm{sph}}$ are smooth.
To this end note that by the linearity of the convolution, the fact that $\mu_i \ast h_j$, $f_{n-1} \ast h_j \in C^{\infty}(\mathbb{S}^{n-1})$, and (\ref{zonalconv}), we have
\[\varphi^j(K) = c_0 + \sum_{i=1}^{n-2}\! \int_{\mathbb{S}^{n-1}}\!\!\!\!\!\! (\mu_i \ast h_j)(u)\, dS_i(K,u) + \int_{\mathbb{S}^{n-1}}\!\!\!\!\!\! (f_{n-1} \ast h_j)(u)\,dS_{n-1}(K,u) + c_n V_n(K).  \]
Thus, an application of Theorem \ref{smoothsphiso} completes the proof. \hfill $\blacksquare$

\vspace{1cm}

\centerline{\large{\bf{ \setcounter{abschnitt}{6}
\arabic{abschnitt}. Integral Geometry of Minkowski Valuations}}}

\reseteqn \alpheqn \setcounter{theorem}{0}

\vspace{0.6cm}

In this final section we apply Theorem \ref{main1} to establish a Crofton formula
for continuous, translation invariant, and $\mathrm{SO}(n)$ equivariant Minkowski valuations. Combining this with Hadwiger's general kinematic formula,
allows us to also deduce a kinematic formula for such Minkowski valuations.

We begin by recalling the classical \emph{Crofton formula} (see, e.g., \textbf{\cite[\textnormal{p.\ 124}]{Klain:Rota}}) for intrinsic volumes: For $0 \leq i, j \leq n$ and $K \in \mathcal{K}^n$, we have
\begin{equation} \label{crofform}
\int_{\mathrm{AGr}_{n-i,n}}\!\!\! V_j(K \cap E)\,d\sigma_{n-i}(E) = \left[\begin{array}{c} i + j\\ j \end{array}\right] V_{i+j}(K).
\end{equation}
Here, $\mathrm{AGr}_{i,n}$ denotes the affine Grassmannian of $i$ planes in $\mathbb{R}^n$ and $\sigma_{i}$
is the rigid motion invariant measure on $\mathrm{AGr}_{i,n}$ normalized such that the set of planes having non-empty intersection with the Euclidean unit ball in $\mathbb{R}^n$ has measure
\[\left[\begin{array}{c} n\\i \end{array}\right]\kappa_{n-i} := {n \choose i} \frac{\kappa_n}{\kappa_i}.  \]

The Crofton formula (\ref{crofform}) is intimately related with the \emph{general kinematic formula}: For $0 \leq j \leq n$ and $K, L \in \mathcal{K}^n$, we have
\begin{equation} \label{genkinform}
\int_{\overline{\mathrm{SO}(n)}} \!\!\! V_j(K \cap gL)\,dg = \sum_{i=0}^{n-j} \left[\begin{array}{c} i + j\\ j \end{array}\right] \left[\begin{array}{c} n \\ i \end{array}\right]^{-1} V_{i+j}(K)V_{n-i}(L),
\end{equation}
where $\overline{\mathrm{SO}(n)} = \mathrm{SO}(n) \ltimes \mathbb{R}^n$ (see the books \textbf{\cite{Klain:Rota, schneiderweil}} for more information).

The obvious connection between (\ref{crofform}) and (\ref{genkinform}) is just a special case of \emph{Hadwiger's general integral geometric theorem}
(see, \textbf{\cite[\textnormal{p.\ 173}]{schneiderweil}}), which states in the translation invariant case that for every $\phi \in \mathbf{Val}$ and $K, L \in \mathcal{K}^n$, we have
\begin{equation} \label{hadgenintgeothm}
\int_{\overline{\mathrm{SO}(n)}}\!\!\! \phi(K \cap gL)\,dg = \sum_{i=0}^n V_{n-i}(L) \left[\begin{array}{c} n \\ i \end{array}\right]^{-1}\!\!\int_{\mathrm{AGr}_{n-i,n}}\!\!\! \phi(K \cap E)\,d\sigma_{n-i}(E) .
\end{equation}

An application of (\ref{hadgenintgeothm}) to the real valued associated valuation of an $\mathrm{SO}(n)$ equivariant Minkowski valuation $\Phi \in \mathbf{MVal}$
immediately yields the following kinematic formula for such Minkowski valuations.

\begin{koro} If $\Phi \in \mathbf{MVal}$ is $\mathrm{SO}(n)$ equivariant, then
\begin{equation} \label{kinformminkval}
\int_{\overline{\mathrm{SO}(n)}}\!\!\! h_{\Phi(K \cap gL)}(u)\,dg = \sum_{i=0}^n V_{n-i}(L) \left[\begin{array}{c} n \\ i \end{array}\right]^{-1}\!\!\int_{\mathrm{AGr}_{n-i,n}}\!\!\!\!\!\!\!\!\!
h_{\Phi(K \cap E)}(u)\,d\sigma_{n-i}(E) \end{equation}
for every $K, L \in \mathcal{K}^n$ and $u \in \mathbb{S}^{n-1}$.
\end{koro}

Note that the sum on the right hand side of (\ref{kinformminkval}) is again the support function of a convex body.
Thus, it remains to determine the Crofton integral in (\ref{kinformminkval}). In view of Lemma \ref{lemdec} and Theorem \ref{main1}, this is accomplished by our final result.

\begin{theorem} \label{thmcrofminkval} Suppose that $1 \leq j \leq n - 2$ and $1 \leq i \leq n - j - 1$. If $F_j \in \mathbf{CVal}_i$ is
$\mathrm{SO}(n)$ equivariant and, for $K \in \mathcal{K}^n$, given by
\[F_{j,K} = S_j(K,\cdot) \ast \mu   \]
for some (uniquely determined) $\mathrm{SO}(n - 1)$ invariant measure $\mu \in \mathcal{M}_{\mathrm{o}}(\mathbb{S}^{n-1})$, then
\begin{equation} \label{laststate17}
\int_{\mathrm{AGr}_{n-i,n}}\!\!\!\!\!\!\!\!\! F_{j,K \cap E}\,d\sigma_{n-i}(E) = q_{n,i,j}\,S_{i+j}(K,\cdot) \ast (\mu \ast \Box_{n-j+1}\breve{g}_{n-i-j+1}),
\end{equation}
where $q_{n,i,j} = \frac{2^i}{i!\kappa_i}\prod_{k=j}^{i+j-1}c_{n,k}$ with
\[c_{n,k} = \frac{k(n-k-1)(n-k+1)\kappa_{n-k-2}^2 \kappa_{n-k+1}\kappa_k}{2(n-k)(k+1)\kappa_{n-k-3}\kappa_{n-k}^2 \kappa_{k-1}}.  \]
\end{theorem}

\noindent {\it Proof.} Consider the isomorphism $\Theta_j: C_{\mathrm{o}}^{\infty}(\mathbb{S}^{n-1}) \rightarrow C_{\mathrm{o}}^{\infty}(\mathbb{S}^{n-1})$, defined by
\[\Theta_jf = c_{n,j\,} \Box_{n-j+1} f \ast \breve{g}_{n-j} = c_{n,j\,} f \ast \Box_{n-j+1} \breve{g}_{n-j}.  \]
Here and in (\ref{laststate17}), $\Box_{k}\breve{g}_{l}$ is to be understood in the sense of distributions, where we use the canonical extension of the selfadjoint operator
$\Box_{k}$ to $C_{\mathrm{o}}^{-\infty}(\mathbb{S}^{n-1})$.

Now, let us first assume that $F_j$ is smooth, that is, $\mu$ is absolutely continuous with respect to spherical Lebesgue measure with density $f \in C_{\mathrm{o}}^{\infty}(\mathbb{S}^{n-1})$.
In this case it was proved by the authors in \textbf{\cite[\textnormal{Theorem 6.3}]{SchuWan13}} that
\begin{equation} \label{Lopsmooth}
 \int_{\mathrm{AGr}_{n-1,n}}\!\!\!\!\!\!\!\!\! F_{j,K \cap E}\,d\sigma_{n-1}(E) = S_{j+1}(K,\cdot) \ast \Theta_j f.
\end{equation}
In order to obtain from this the more general formula (\ref{laststate17}), we use the following well known relation (which can be proved by induction using Crofton's formula;
see, e.g., \textbf{\cite[\textnormal{p.\ 124}]{Klain:Rota}})
\[ \int_{\mathrm{AGr}_{n-i,n}} \!\!\!\!\!\!\!\!\!\!\!\!\!\!\! f(E)\,d\sigma_{n-i}(E) = \frac{2^i}{i!\kappa_i}\!
\int_{\mathrm{AGr}_{n-1,n}} \!\!\!\!\!\!\!\! \cdots \int_{\mathrm{AGr}_{n-1,n}}\!\!\!\!\!\!\!\!\!\!\!\!\!\!\! f(E_1 \cap \cdots \cap E_i)\,d\sigma_{n-1}(E_1) \cdots d\sigma_{n-1}(E_i)    \]
for every Borel measurable $f \in L^1(\mathrm{AGr}_{n-i,n})$. Comparing this with (\ref{Lopsmooth}), we obtain
\begin{equation} \label{laststate18}
\int_{\mathrm{AGr}_{n-i,n}}\!\!\!\!\!\!\!\!\! F_{j,K \cap E}\,d\sigma_{n-i}(E) = \frac{2^i}{i!\kappa_i} S_{i+j}(K,\cdot) \ast \Theta_{i+j-1}\cdots \Theta_{j+1}\Theta_j f.
\end{equation}
Next, note that if $\tau_{\bar{e}} = \delta_{\bar{e}}- \pi_1\delta_{\bar{e}} \in \mathcal{M}_{\mathrm{o}}(\mathbb{S}^{n-1})$, where $\delta_{\bar{e}}$ is
the Dirac measure supported in $\bar{e} \in \mathbb{S}^{n-1}$, then, by (\ref{funkheckgen}), $f \ast \tau_{\bar{e}} = f$ for every
$f \in C_{\mathrm{o}}^{\infty}(\mathbb{S}^{n-1})$. But since $\Box_k \breve{g}_k = \tau_{\bar{e}}$, (\ref{laststate17}) follows from (\ref{laststate18}) and the definition of $\Theta_j$.

In order to establish (\ref{laststate17}) in the general case, where $F_j$ is merely continuous, we use a spherical approximate identity $h_k$, $k \in \mathbb{N}$, (instead of repeating the arguments from the proof
of \textbf{\cite[\textnormal{Theorem 6.3}]{SchuWan13}}) to define $F_{j,K}^k = F_{j,K} \ast h_k$ for every $K \in \mathcal{K}^n$. Then, $F_j^k \in \mathbf{CVal}_j$ is
$\mathrm{SO}(n)$ equivariant and smooth and, by what we have already shown and the fact that multiplier transformations commute,
\[\int_{\mathrm{AGr}_{n-i,n}}\!\!\!\!\!\!\!\!\! F^k_{j,K \cap E}\,d\sigma_{n-i}(E) = q_{n,i,j}\,S_{i+j}(K,\cdot) \ast (\mu \ast \Box_{n-j+1}\breve{g}_{n-i-j+1}) \ast h_k.\]
Letting now $k \rightarrow \infty$, we obtain (\ref{laststate17}) from Lemmas \ref{weakcontconv} and \ref{approxid}. \hfill $\blacksquare$

\vspace{0.3cm}

We conclude with the remark that equivalent forms of Theorem \ref{thmcrofminkval} were obtained very recently, independently, and using different approaches by
Bernig and Hug \textbf{\cite{bernighug2015}} and Goodey, Hug, and Weil \textbf{\cite{goodeyhugweil15}}.

\vspace{0.5cm}

\noindent {{\bf Acknowledgments} The first author was
supported by the European Research Council (ERC), Project number: 306445, and the Austrian Science Fund (FWF), Project number:
Y603-N26. The second author was supported by the German Research Foundation (DFG), Project number: BE 2484/5-1.

\begin{small}

\[ \begin{array}{ll} \mbox{Franz Schuster} & \mbox{Thomas Wannerer} \\
\mbox{Vienna University of Technology \phantom{wwwwWW}} & \mbox{Goethe University Frankfurt} \\ \mbox{franz.schuster@tuwien.ac.at} & \mbox{wannerer@mathematik.uni-frankfurt.de}
\end{array}\]

\end{small}

\end{document}